\documentclass[12pt]{article}
\usepackage{amsmath}
\usepackage{enumerate}
\usepackage{amsthm}
\usepackage{authblk}
\usepackage{adjustbox}
\usepackage{graphicx}
\usepackage{blkarray}
\usepackage[backend=biber,style=apa,citestyle=authoryear]{biblatex}
\usepackage{amssymb}
\usepackage{float}
\usepackage{commath} 
\usepackage{mathrsfs}
\usepackage{csquotes}
\usepackage{tikz-cd}
\usepackage{algorithm}
\usepackage{algpseudocode}
\usepackage{subcaption}
\usepackage{mathtools}
\usepackage[english]{babel}

\newtheorem{theorem}{Theorem}[section]
\newtheorem{lemma}[theorem]{Lemma}
\newtheorem{prop}{Proposition}[section]

\newcommand{\blind}{1}

\begin{document}
\def\spacingset#1{\renewcommand{\baselinestretch}%
{#1}\small\normalsize} \spacingset{1}


\if1\blind
{
  \title{\bf Diffusion on the circle and a stochastic correlation model}
  \author{Sourav Majumdar\\
    Department of Management Sciences, Indian Institute of Technology Kanpur\\
    souravm@iitk.ac.in\\
    and \\
    Arnab Kumar Laha \\
    Operations and Decision Sciences, Indian Institute of Management Ahmedabad\\
    arnab@iima.ac.in}
  \maketitle

} \fi

\if0\blind
{
  \bigskip
  \bigskip
  \bigskip
  \begin{center}
    {\LARGE\bf Diffusion on the circle and a stochastic correlation model}
\end{center}
  \medskip
} \fi

\bigskip
\begin{abstract}
We develop diffusion models for time-varying correlation using stochastic processes defined on the unit circle. Specifically, we study Brownian motion on the circle and the von Mises diffusion, and propose their use as continuous-time models for correlation dynamics. The von Mises process, introduced by Kent (1975) as a characterization of the von Mises distribution in circular statistics, does not have a known closed-form transition density, which has limited its use in likelihood-based inference. We derive an accurate analytical approximation to the transition density of the von Mises diffusion, enabling practical likelihood-based estimation. We study inference for discretely observed circular diffusions, establish consistency and asymptotic normality of the resulting estimators, and propose a stochastic correlation model for financial applications. The methodology is illustrated through simulation studies and empirical applications to equity-foreign exchange market data.\end{abstract}

\noindent%
{\it Keywords:} Directional Statistics, Statistical inference of continuous-time processes, von Mises process, Quantitative Finance
\vfill

\newpage
\spacingset{1.9} 

\section{Introduction}\label{sec: intro}

The von Mises distribution is a fundamental probability model for circular data, that is, for observations taking values on the unit circle. It was introduced in 1918 by Richard von Mises in his study of the so-called ``integralness" of atomic weights (\cite{vonMises1918Ganzzahligkeit}). Motivated by the question of whether observed atomic weights could be regarded as integer multiples of a fundamental unit, von Mises reformulated the problem by mapping the fractional parts of the measurements onto the unit circle. He then sought a probability distribution for circular errors with the defining property that the maximum likelihood estimator of the location parameter coincides with the defined mean direction. This requirement uniquely leads to what is now known as the von Mises distribution, which is widely regarded as the circular analogue of the normal distribution on the real line.

Since its introduction, several alternative characterizations of the von Mises distribution have been proposed, including maximum entropy arguments and conditioning constructions; see, for example, \cite{mardia2000directional} for a comprehensive overview.

One important way in which the von Mises distribution arises is as the stationary distribution of a mean-reverting diffusion on the circle. The von Mises process, defined by the stochastic differential equation,

$$
d\theta_t=-\lambda\sin(\theta_t-\mu)dt+\sigma dW_t
$$

was introduced by Kent (\cite{kent1975disc,kent1978time}) as a stochastic process characterization of the von Mises distribution. As a continuous-time Markov process, the transition density of this diffusion fully determines its probabilistic behaviour and is essential for statistical inference and applied modelling. However, despite its importance, a closed-form expression for the transition density of the von Mises process is not available. Recent works (\cite{garcia2019langevin,GarcaPortugus2025}) explicitly note that the transition density remains unknown.

In this article, we derive an accurate approximation to the transition density of the von Mises process. This approximation facilitates likelihood-based inference and enables the practical use of the von Mises diffusion in applied settings. In particular, it broadens the applicability of circular diffusion models in areas such as meteorology, where wind direction is naturally modelled on the circle. We also explore a novel application of the von Mises process in financial modelling.

A pertinent problem in quantitative finance is the joint modelling of multiple assets. Such issues arise in empirical analysis of financial data, portfolio optimization, and multi-asset derivatives pricing, among others. The usual modelling approach is to assume that the dependence between the assets is constant over time. However, empirical evidence suggests that the correlation in stock markets is stochastic (\cite{ball2000stochastic,driessen2009price,faria2022correlation}).
This has led to significant interest towards the study of estimation of quadratic co-variance and its time varying nature, see \cite{ait2010high, fan2016incorporating} and the references therein. One is also interested in parametric models for continuous-time correlation to perform no-arbitrage pricing of derivatives and forecasting. One of the first proposals to model continuous-time correlation was considered in \cite{vanemmerich2006modelling}, wherein a bounded functional of a Brownian motion was used to model correlation. Several extensions based on this approach have been considered \cite{teng2015pricing,teng2016modelling,teng2016heston}. \cite{teng2018pricing} consider a Jacobi process to model correlation. Another method of modeling covariance is via Wishart processes, which has also received a lot of attention in the literature \cite{fonseca2007option,gourieroux2010derivative,philipov2006multivariate,asai2006multivariate}.

The transition density of a diffusion is not always available in closed-form, which makes parametric inference of the diffusion difficult. This has led to significant interest towards the study of numerical and simulation based approaches to estimating the parameters of a diffusion, see \cite{kou2012multiresolution,elerian2001likelihood,ait2002maximum,beskos2006exact}. See \cite{kalogeropoulos2011likelihood} for a simulation method for estimating correlated diffusions. For a recent survey, see \cite{craigmile2023statistical}.

Circular statistics studies data sampled from the unit circle,
\cite{jammalamadaka2001topics,mardia2000directional}.
In this article, we study the estimation of some
diffusions on the unit circle and explore a novel
application to quantitative finance in stochastic
correlation modeling. We consider two circular
diffusions in this paper- the Circular Brownian Motion
and the von Mises Process.
We are interested in estimating the parameters of a discretely observed circular diffusion through maximum likelihood estimation. The analytical form of the transition density of the von Mises process allows us to compute MLEs without relying on fully numerical or simulation approaches for the estimation of parameters of SDEs. Recently \cite{garcia2019langevin} also noted the lack of knowledge of the transition density of the von Mises process in the literature. They proposed numerical approaches to estimate the von Mises process.

In this article, we consider a new approach to model correlation using circular diffusions. If $X$ and $Y$ are two random vectors in $\mathbb{R}^n$, then the Pearson correlation coefficient between $X$ and $Y$ is $\rho=\frac{\langle X,Y \rangle}{||X||||Y||}=\cos \theta$, where, $X=(x_1,x_2,\ldots,x_n),Y=(y_1,y_2,\ldots,y_n)$, $\langle X,Y \rangle=\sum_{i=1}^nx_iy_i$ and $||X||=\sqrt{x_1^2+x_2^2+\ldots +x_n^2}$. This leads to a natural choice to model $\rho$ by modelling $\theta$ as a process on the circle. While ideas from circular statistics have been applied earlier to the statistical analysis of financial data (see \cite{sengupta2019universal,sengupta2023circular}), the use of circular diffusion as a model for correlation has not been considered to the best of our knowledge.

The rest of the article is structured as follows. In Section \ref{sec: circdiff} we discuss some background on diffusions on the circle. We derive an approximate transition density of the von Mises process and investigate its accuracy. In Section \ref{sec: estimation}, we discuss the estimation of diffusions on the circle. We perform a simulation study and also perform real data analysis using them. In Section \ref{sec: stochcorr} we propose a stochastic correlation model, we discuss its estimation and a bootstrap procedure. In Section \ref{sec: application} we fit the proposed stochastic correlation model to real financial data and discuss the analysis. In Section \ref{sec: conclusion} we conclude the article.

\section{Circular Diffusion}\label{sec: circdiff}

We first collect some known facts about Brownian motion on the circle. Let $\mathbb{S}^1$ denote the circle group, which is also a Lie group (see \cite{Stillwell_2008} for a background on Lie groups). \cite{ito1950brownian} studied Brownian motion on Lie groups. Let $\theta$ denote the coordinate. For a Brownian motion on Lie group, it is shown that its transition density is uniquely characterized by its infinitestimal generator. $p(\theta_u;\theta_v)$ denotes the transition density function, i.e. the probability density of $\theta_u$ given $\theta_v$, where $\theta_u$ is the position of the particle at time $u$. The infinitesimal generator (see \cite{ito1950brownian}, (1.2)) for a Brownian motion on the circle group is given by,
\begin{equation}
    D=\frac{1}{2}\sigma^2\frac{\partial^2}{\partial \theta^2}
\end{equation}
and its forward equation (see \cite{ito1950brownian}, Theorem 2) is,
\begin{equation}\label{eq:cbm_forward}
    \frac{\partial}{\partial t}p(\theta_t;\theta_0)=Dp(\theta_t;\theta_0)
\end{equation}
See \cite{oksendal2013stochastic}, p. 117 for a definition of the infinitesimal generator. It is further shown that the paths of the Brownian motion on circle is continuous in the sense of Kolmogorov-Feller (\cite{ito1950brownian}) and its SDE representation is given by (\cite{ito1950brownian}, (2.6)),
\begin{equation}
    d\theta_t=\sigma dB_t
\end{equation}

where $\theta_t$ denotes the position of the Brownian motion on the circle at time $t$. 

One can also characterize the Brownian motion on the circle alternatively, see \cite{mardia2000directional}, p. 51. Let $B_t$, $B_0=0$ denote the standard Brownian motion on $\mathbb{R}$. Then $\theta_t \in [0,2\pi)$ is said to follow a circular Brownian motion as the process obtained by wrapping the Brownian motion on the circle, i.e., $\theta_t=(\sigma B_t \text{ mod } 2\pi)$.

Equivalently it can also be characterized as following where $\theta_t$ is said to the follow the circular Brownian motion if it satisfies the following SDE,
\begin{equation}
	d\theta_t=\sigma dB_t
\end{equation}

with periodic boundary conditions, $p(\theta_u;\theta_v)=p(\theta_u+2k\pi;\theta_v)$, where $k \in \mathbb{Z},u>v>0$. The initial condition is $p_0(\theta)=\delta(\theta-0)$, where $\delta$ denotes the Dirac-Delta function. The transition density is given in (\ref{cbm_tpd}), by the Kolmogorov-Feller criterion the diffusion is continuous.

$\theta_t$ is said to follow the von Mises process (\cite{kent1975disc}) if it satisfies the following SDE,
 \begin{equation}
	d\theta_t=-\lambda\sin(\theta_t-\mu)dt+\sigma dB_t
\end{equation}

for $\lambda>0,\sigma>0$ and periodic boundary conditions.

The von Mises process was proposed as the circular analogue of the Ornstein-Uhlenbeck process on $\mathbb{R}$ in \cite{kent1975disc}. It is mean-reverting in the sense that the drift term $-\lambda\sin(\theta_t-\mu)$ pulls the particle $\mathbb{\theta}_t$ towards the mean direction $\mu$. It can be shown that the stationary distribution of the von Mises process is the von Mises distribution with density $f(\theta)$,
\begin{equation}
	f(\theta)=\frac{1}{2\pi I_0(\kappa)}\exp(\kappa\cos(\theta-\mu)), 0\leq \theta < 2\pi
\end{equation}

where $\kappa=\frac{2\lambda}{\sigma^2}$ and $I_0(\cdot)$ is the zeroth-order modified Bessel function of the first kind. The von Mises distribution is widely used in circular statistics. To the best of our knowledge, no analytical solution of the transition density of the von Mises process has been reported in the literature. In this article, we present an approximate analytical solution for the transition density of the von Mises process. 

\cite{kent1978time} presented a construction through which one can obtain a time-reversible diffusion with a specified stationary distribution. We collect a special case of this construction in the following Lemma,

\begin{lemma}\label{kents_construction}
    Let $f(\theta)$ be a strictly positive density on the circle, then the diffusion supported on the circle whose stationary density is $f(\theta)$ is given by,
    \begin{equation}
    d\theta_t=\frac{\sigma^2}{2}\frac{\partial}{\partial \theta_t} \log f(\theta_t)dt+\sigma dB_t
    \end{equation}
    satisfying periodic boundary conditions.
\end{lemma}
\begin{proof}
	See the appendix for a proof.
\end{proof}
Applying Lemma \ref{kents_construction}, we construct the circular Brownian motion and the von Mises process.

For $f(x)=\frac{1}{2\pi}$, the density of the uniform distribution on the circle; applying Lemma \ref{kents_construction} and  the corresponding SDE given by,

\begin{equation}
	d\theta_t=\sigma dB_t
\end{equation}

alongwith the periodic boundary condition is referred to as the circular Brownian motion.

The forward equation for the circular Brownian motion is given by, 

\begin{equation}
\frac{\partial}{\partial t} p(\theta_t;\theta_0)=\frac{\sigma^2}{2}\frac{\partial^2}{\partial \theta^2}
 p(\theta_t;\theta_0)\end{equation}

It can be checked by solving the forward equation of the circular Brownian motion (see (\ref{eq:cbm_forward})) that the transition density of $\theta_t$, $0<s<t$, is,
\begin{equation}\label{cbm_tpd}
    p(\theta_t;\theta_s)=\frac{1}{\sigma\sqrt{2\pi(t-s)}}\sum_{k=-\infty}^{\infty}\exp\left(\frac{-(\theta_t-\theta_s+2k\pi)^2}{2\sigma^2(t-s)}\right).
\end{equation}
which is the well-known wrapped normal distribution on the circle with parameters $\theta_s$ and  $\sigma\sqrt{t-s}$, denoted by $\text{WN}(\theta_s,\sigma\sqrt{t-s})$. Thus, the stationary distribution of the circular Brownian motion is the uniform distribution and the transition density is the wrapped normal distribution.

For $f(x)=\frac{1}{2\pi I_0(\kappa)}\exp(\kappa\cos(\theta-\mu))$, the density of the von Mises distribution with mean direction $\mu$ and concentration  $\kappa=\frac{2\lambda}{\sigma^2}$, where $\lambda$ is the mean-reversion parameter given in (\ref{sde: vmp}); applying Lemma \ref{kents_construction} and the corresponding SDE is given by,
\begin{equation}\label{sde: vmp}
	d\theta_t=-\lambda\sin(\theta_t-\mu)dt+\sigma dB_t
\end{equation}

where $\theta,\mu\in [0,2\pi), \lambda,\sigma>0$. This SDE is also referred to as the von Mises process.

The corresponding Kolmogorov-Forward equation is given by,
\begin{equation}\label{vm_forward}
\frac{\partial}{\partial t} p(\theta_t;\theta_0)=\frac{\sigma^2}{2}\frac{\partial^2}{\partial \theta^2} p(\theta_t;\theta_0)+\lambda\sin(\theta_t-\mu)\frac{\partial}{\partial \theta}p(\theta_t;\theta_0)+\lambda\cos(\theta_t-\mu)p(\theta_t;\theta_0)
\end{equation}

The stationary distribution of the von Mises process is the von Mises distribution as noted in its construction. However, the transition density of the von Mises process has not yet been reported in the literature.

\subsection{Approximate transition density of von Mises process}\label{subsec: approx_tpd}

We use the approximation framework presented in \cite{martin2019analytical} to obtain an approximate analytical solution for the forward equation of the von Mises process.

For deriving the approximate transition density of the von Mises process, we use the fact that the transition density of the Ornstein-Uhlenbeck process is known in closed form. A logarithmic transformation of the forward equation is considered, and the Ornstein-Uhlenbeck density then suggests an ansatz for more general forms of the drift term.

Therefore we state the following result,

\begin{theorem}\label{thm:vmp}
\begin{align}\label{tpd_vm}
	p(\theta_t;\theta_0)&\propto\sum_{k=-\infty}^{\infty}\exp\left(\frac{-\gamma\sqrt{q}(\theta_t+2\pi k-\theta_0)^2}{2(1-q)}\right)\nonumber\\
			    &\left(\frac{1}{2\pi I_0(\kappa)}\right)^{\frac{1-\sqrt{q}}{1+\sqrt{q}}}\exp\left(\frac{\kappa(\cos(\theta_t-\mu)-\sqrt{q}\cos(\theta_0-\mu))}{1+\sqrt{q}}\right)	
\end{align}

where $\gamma=\kappa(I_1(\kappa)/I_0(\kappa))$ and $q=\exp(-\gamma\sigma^2t)$. $p(\theta_t;\theta_0)$ is the approximate transition density of the von Mises process, in the sense that for $t\to 0$, the error term is $o(1)$ and for $t\to\infty$, $p(\theta;\theta_0,t)\to f_e$. Here $f_e$ is the von Mises distribution.
\end{theorem}
\begin{proof}
	See the appendix for a proof.
\end{proof}
We note that as $t\to \infty$,  $p(\theta_t;\theta_0)\to f_e$, which is the expected theoretical stationary distribution. The constant of proportionality $(C(\theta_0,\lambda,\sigma,\mu,t))$ can be obtained by numerically integrating $p(\theta;\theta_0,t)$ over  $\theta$, i.e.,

\begin{equation}
	C(\theta_0,\lambda,\sigma,\mu,t)=\frac{1}{\int p(\theta_t;\theta_0)d\theta}
\end{equation}

\subsection{Numerical experiment}

In (\ref{tpd_vm}) we obtain an approximate analytical solution of the von Mises process. To validate the accuracy of the approximate analytical solution we compare it with a numerical solution of the transition density.

One possible approach is to solve the forward equation (\ref{vm_forward}) numerically. Towards this a standard approach is the Crank-Nicholson scheme which was also considered in Section 3.4.1, \cite{garcia2019langevin}. In this approach a uniform discretization of the space, i.e., $[0,2\pi)$ is considered into $k$ points and similarly time $(0,T]$ is discretized uniformly into  $m$ points. The PDE is also further discretized and the resulting numerical estimate of the density is obtained at these  $k$ points for each  $m$ time points. The reader may refer to Section 3.4.1, \cite{garcia2019langevin} for an exposition of the Crank-Nicholson scheme.  

Let $p$ and $q$ denote two probability density functions, then the Hellinger distance (see \cite{le2000asymptotics}) is defined as,

\begin{equation}
	h(p,q)=\frac{1}{\sqrt{2}}\sqrt{\int_{\mathbb{R}}(\sqrt{dp}-\sqrt{dq})^2}
\end{equation}

Hellinger distance is used to quantify the similarity between two probability densities. It is known that $0\leq h(p,q) \leq 1$, where if $h(p,q)$ is close to 0, then  $p$ and  $q$ are interpreted to be similar.

Let $p$ denote the density of the approximate analytical solution and let $q$ denote the density of the numerical solution. Since the numerical solution $q$ via the Crank-Nicholson scheme is discrete, we consider the discrete version of the Hellinger distance by discretizing $p$ at the same $k$ points. This is given by,
\begin{equation}
    h(p,q)=\frac{1}{\sqrt{2}}\sqrt{\sum_{i=1}^k(\sqrt{p_i}-\sqrt{q_i})^2}
\end{equation}

\begin{figure}[H]
\centering
\begin{adjustbox}{minipage=1.2\textwidth,scale=0.6}
\begin{subfigure}{\textwidth}
\includegraphics[width=\textwidth]{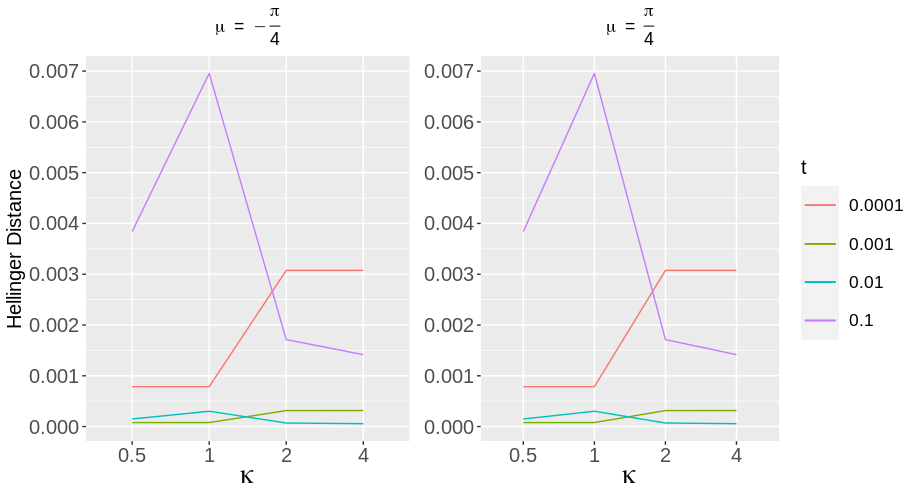}
\caption{}
\label{subfig:hell_a}
\end{subfigure}
\vfill
\begin{subfigure}{\textwidth}
\includegraphics[width=\textwidth]{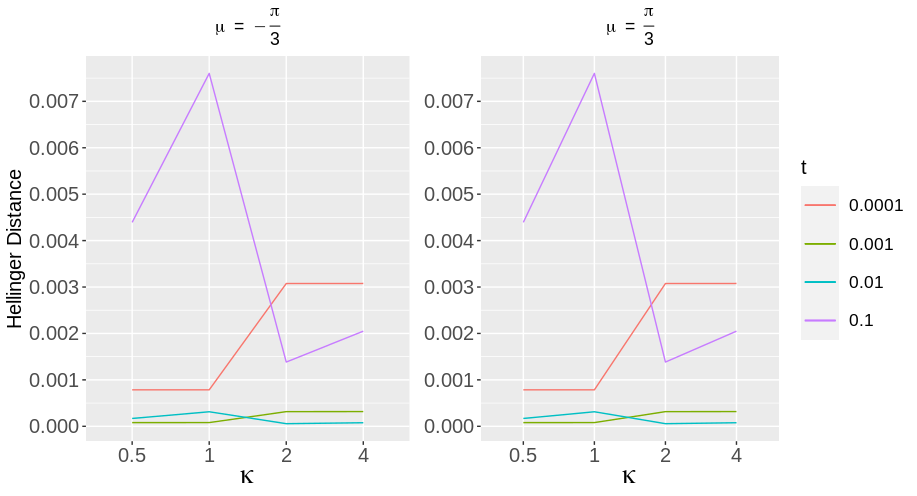}
\caption{}
\label{subfig:hell_b}
\end{subfigure}
\vfill
\begin{subfigure}{\textwidth}
\includegraphics[width=\textwidth]{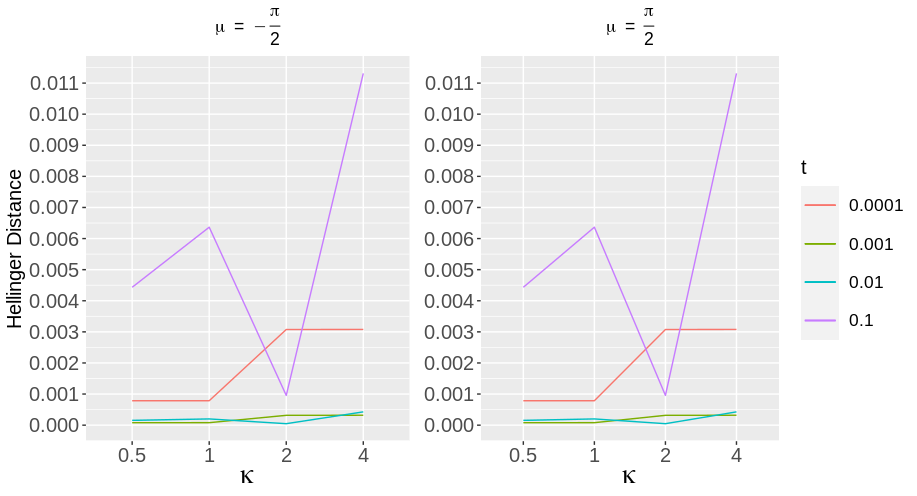}
\caption{}
\label{subfig:hell_c}
\end{subfigure}
\caption{Hellinger distance between the analytical $(p(\theta;\theta_0,t))$ and numerical approximation of the forward equation of the von Mises process. The parameters are $\theta_0=0$ and  $\kappa=0.5$ $(\lambda=1,\sigma=2)$,  $\kappa=1$ $(\lambda=2,\sigma=1)$,  $\kappa=2$  $(\lambda=1,\sigma=2)$ and  $\kappa=4$  $(\lambda=2,\sigma=2)$. We compare the Hellinger distance at $\mu=-\frac{\pi}{2},-\frac{\pi}{3},-\frac{\pi}{4},\frac{\pi}{4},\frac{\pi}{3},\frac{\pi}{2}$. For a given $\mu$ we plot the Hellinger distance between the numerical and the analytical density for different values of  $t$, where  $t$ is the time elapsed. For the Crank-Nicholson scheme for the numerical solution we set $k=3000,m=20000$.}
\label{fig:hell}
\end{adjustbox}
\end{figure}

We use the function \texttt{dTpdPde1D} in the R package \texttt{sdetorus} (\cite{sdetorus}) to utilize the implementation of Crank-Nicholson scheme to numerically solve the forward equation of the von Mises process. We compare analytical approximate solution $p(\theta;\theta_0,t)$ (\ref{tpd_vm}) with the numerical solution via the Hellinger distance between them. We set $k=3000,m=20000$ for the Crank-Nicholson scheme. We set  $\theta_0=0$ and we compare the densities after time $t$ has elapsed, where  $t\in\{0.0001,0.001,0.01,0.1\}$ for varying $\mu\in\{-\frac{\pi}{2},-\frac{\pi}{3},-\frac{\pi}{4},\frac{\pi}{4},\frac{\pi}{3},\frac{\pi}{2}\}$. For each $t$ and  $\mu$ we also evaluate the Hellinger distance for different values of  $\kappa$, i.e., $\kappa=0.5$ $(\lambda=1,\sigma=2)$,  $\kappa=1$ $(\lambda=2,\sigma=1)$,  $\kappa=2$  $(\lambda=1,\sigma=1)$ and  $\kappa=4$  $(\lambda=2,\sigma=2)$. We report the results in Figure \ref{fig:hell}. We find that for all cases considered the Hellinger distance is very low $(\lessapprox 0.01)$, which indicates the accuracy of the approximate analytical solution.

\section{Estimation of Circular Diffusions}\label{sec: estimation}

In this section we consider the estimation of a discretely observed Circular diffusion. Let $\theta_t$ be a circular diffusion and we observe $\theta_t$ at times  $0=t_1<t_2<\cdots<t_n=T$. We consider the case where $\theta_t$ is the circular Brownian motion and the von Mises process.

The quadratic variation of $\theta_t$, where $\theta_t$ could be either the circular Brownian motion or the von Mises process, is given by,

 \begin{equation}
	 [\theta]_t=\int_0^T\sigma^2dt
\end{equation}

Hence, $\sigma^2=\frac{[\theta_t]}{T}$. The sample quadratic variation estimator is given by,

\begin{equation}
	[\theta]_t^*=\sum_{i=1}^{n-1}(\theta_i^*)^2
\end{equation}

where $\theta_i^*=(\theta_{t_{i+1}}-\theta_{t_i}) \text{ mod } \pi$. It is known that $[\theta]_t^*\xrightarrow{ucp}[\theta]_t$ as $\sup_i|t_{i+1}-t_{i}|\to 0$, see \cite{protter2005}, Theorem 22, p. 66. Therefore the diffusion coefficient estimator is,

\begin{equation}\label{quad_var_est}
	\hat{\sigma}^2=\frac{[\theta]_t^*}{T}
\end{equation}

If $\sup_i|t_{i+1}-t_i|$ is sufficiently small, i.e.\ we sample from the diffusion very frequently, then the estimator is expected to perform well.

We note that the consistency and asymptotic normality of the maximum likelihood estimation (MLE) of a discretely observed It{\^o} process has been established, see~\cite{prakasa1983asymptotic,yoshida1992estimation}.  To estimate the parameters in the drift term of the von Mises process, i.e. $\lambda,\mu$, we propose to perform numerical estimates of the parameters. We are unable to obtain analytical solutions for $\lambda,\mu$ since we know (\ref{tpd_vm}) only upto proportionality. A promising approach to estimating parameters of non-normalized models is the score matching estimator, see~\cite{hyvarinen2005estimation, mardia2016score}. The score-matching formulation for this model was tried but proved to be not analytically tractable for us.

The log-likelihood of the von-Mises process is given by,

\begin{equation*}
L(\lambda,\mu;\boldsymbol{\theta})
=
\sum_{i=1}^{n-1}
\left[
\log C_i
+
\log p\left(\theta_{i+1}\mid \theta_i;\Delta t\right)
\right],
\end{equation*}

where $C_i$ is a shorthand for  $C(\theta_i,\lambda,\sigma,\mu,t_i)$ and we use  $\hat{\sigma}$ as a plugin estimator for the diffusion coefficient. Then $\arg\max_{\lambda,\mu} L$ to obtain the estimates for  $\lambda,\mu$. 

\begin{theorem}\label{thm:consistency}
Consider the von Mises process
\begin{equation*}
d\theta_t = -\lambda \sin(\theta_t - \mu)\, dt + \sigma\, dW_t, \qquad t \in [0,T],
\end{equation*}
where $\lambda > 0$, $\mu \in [0,2\pi)$, and $\sigma>0$. Let
$\widehat{\lambda},\widehat{\mu},\widehat{\sigma}$ denote the maximum likelihood estimators based on continuous observations on $[0,T]$.

Let $\Theta_\lambda \subset (0,\infty)$, $\Theta_\mu \subset [0,2\pi)$, 
$\Theta_\sigma \subset (0,\infty)$ be open parameter sets, and let 
$K_\lambda \subset \Theta_\lambda$, $K_\mu \subset \Theta_\mu$, 
$K_\sigma \subset \Theta_\sigma$ be compact subsets. Then for every $\varepsilon>0$,
\begin{align*}
\lim_{T\to\infty} \sup_{\lambda\in K_\lambda} 
\mathbb{P}_\lambda\!\left(|\widehat{\lambda}-\lambda|>\varepsilon\right)=0\\
\lim_{T\to\infty} \sup_{\mu\in K_\mu} 
\mathbb{P}_\mu\!\left(|\widehat{\mu}-\mu|>\varepsilon\right)=0\\
\lim_{T\to\infty} \sup_{\sigma\in K_\sigma} 
\mathbb{P}_\sigma\!\left(|\widehat{\sigma}-\sigma|>\varepsilon\right)=0
\end{align*}
Thus $(\widehat{\lambda},\widehat{\mu},\widehat{\sigma})$ is uniformly consistent on compact subsets.
\end{theorem}

\begin{proof}
	See the appendix for a proof.
\end{proof}

\begin{theorem}\label{thm:normality}
Let $(\hat\lambda,\hat\mu)$ denote the maximum likelihood estimators of
$(\lambda,\mu)$ based on continuous observation of the von Mises process,
$$
d\theta_t=-\lambda\sin(\theta_t-\mu)\,dt+\sigma\,dW_t,
\qquad t\in[0,T],
$$
with $\sigma>0$ known.
Then, uniformly on compact subsets of the parameter space,
\begin{align*}
\sqrt{T}(\hat{\lambda}-\lambda)
&\xrightarrow[]{d}
N\!\left(0,\frac{2\sigma^2 I_0(\kappa)}{I_0(\kappa)-I_2(\kappa)}\right),\\[6pt]
\sqrt{T}(\hat{\mu}-\mu)
&\xrightarrow[]{d}
N\!\left(0,\frac{2\sigma^2 I_0(\kappa)}
{\lambda^2\bigl(I_0(\kappa)+I_2(\kappa)\bigr)}\right),
\end{align*}
where $\kappa=2\lambda/\sigma^2$ and $I_\nu$ denotes the modified Bessel
function of the first kind of order $\nu$.
\end{theorem}

\begin{proof}
	See the appendix for a proof.
\end{proof}

\begin{theorem}\label{thm:normality}
Let $\theta=\{\theta_t:0\le t\le T\}$ be continuously observed from the von Mises diffusion
\[
d\theta_t=-\lambda\sin(\theta_t-\mu)\,dt+\sigma\,dW_t,\qquad \theta_t\in\mathbb S^1,
\]
where $\sigma>0$ is known and the parameter $(\lambda,\mu)$ ranges over a compact set
$K\subset (0,\infty)\times[0,2\pi)$.
Let $(\hat\lambda_T,\hat\mu_T)$ be the maximum likelihood estimator of $(\lambda,\mu)$.
Then, uniformly for $(\lambda,\mu)\in K$,
\begin{align*}
\sqrt{T}(\hat{\lambda}_T-\lambda)
&\xrightarrow[]{d}
N\!\left(0,\frac{2\sigma^2 I_0(\kappa)}{I_0(\kappa)-I_2(\kappa)}\right),\\[6pt]
\sqrt{T}(\hat{\mu}_T-\mu)
&\xrightarrow[]{d}
N\!\left(0,\frac{2\sigma^2 I_0(\kappa)}
{\lambda^2\bigl(I_0(\kappa)+I_2(\kappa)\bigr)}\right),
\end{align*}
where $\kappa=2\lambda/\sigma^2$ and $I_\nu$ is the modified Bessel function of the first kind.
\end{theorem}

\subsection{Simulation study}\label{sec:sim_study}

We evaluate the performance of the estimating procedure by simulating from both the circular Brownian motion and the von Mises process. We simulate from the circular Brownian motion and the von Mises process utilizing the Euler-Maruyama discretization of their SDEs and wrapping it. See \cite{Kloeden_1992} for a background on Euler-Maruyama discretization. We vary both the sample size $n$ (i.e., the number of observations) and the time step $\Delta t$ between observations. For each setting, we generate 100 simulated paths and estimate the parameters on each path. In the tables below, we report the empirical bias and the standard deviation of that difference across the 100 replications. 

For the circular Brownian motion we vary $\sigma$ in $\{1,2\}$, $n$ in  $\{1000,5000,10000\}$ and $\Delta t$ in  $\{0.005,0.05,0.5\}$. We report the results in Table \ref{tab:cbm_estimator}. We observe that the Bias of the $\hat{\sigma}$ is close to 0 in general. The variance of  $\sigma-\hat{\sigma}$ decreases as $n$ increases and  $\Delta t$ decreases. For the circular Brownian motion simulations, the quadratic variation estimator for $\sigma$ performs very well, displaying minimal bias and decreasing variance with increasing $n$.

\begin{table}[H]
\scalebox{0.8}{%
\centering
\begin{tabular}{|c|c|c|c|c|} 
\hline
$\sigma$ & $n$     & $\Delta t$    & $\mathbb{E}[\sigma-\hat{\sigma}]$ & $\sqrt{\text{Var}[\sigma-\hat{\sigma}]}$  \\ 
\hline
1     & 1000  & 0.5   & 0.001       & 0.023       \\ 
\hline
1     & 1000  & 0.05  & 0.001       & 0.023       \\ 
\hline
1     & 1000  & 0.005 & 0.000       & 0.021       \\ 
\hline
1     & 5000  & 0.5   & -0.001      & 0.010       \\ 
\hline
1     & 5000  & 0.05  & 0.000       & 0.008       \\ 
\hline
1     & 5000  & 0.005 & 0.003       & 0.010       \\ 
\hline
1     & 10000 & 0.5   & 0.001       & 0.007       \\ 
\hline
1     & 10000 & 0.05  & 0.001       & 0.006       \\ 
\hline
1     & 10000 & 0.005 & 0.000       & 0.007       \\ 
\hline
2     & 1000  & 0.5   & 0.078       & 0.042       \\ 
\hline
2     & 1000  & 0.05  & 0.002       & 0.047       \\ 
\hline
2     & 1000  & 0.005 & -0.009      & 0.047       \\ 
\hline
2     & 5000  & 0.5   & 0.083       & 0.014       \\ 
\hline
2     & 5000  & 0.05  & 0.001       & 0.022       \\ 
\hline
2     & 5000  & 0.005 & -0.003      & 0.019       \\ 
\hline
2     & 10000 & 0.5   & 0.083       & 0.010       \\ 
\hline
2     & 10000 & 0.05  & 0.002       & 0.013       \\ 
\hline
2     & 10000 & 0.005 & -0.001      & 0.015       \\
\hline
\end{tabular}%
}
\caption{Simulation results for the circular Brownian motion. Columns report the empirical mean (bias) and standard deviation (SD) of $\sigma-\hat{\sigma}$ across 100 simulated paths for each combination of $\sigma$, $n$, and $\Delta t$.}
\label{tab:cbm_estimator}
\end{table}

We next consider von Mises process simulations, where the true parameters are $\mu,\lambda$, and $\sigma$. We again vary $n$ in $\{1000, 5000, 10000\}$ and $\Delta t$ in $\{0.005, 0.05, 0.5\}$. We vary the parameter configurations as follows where, $\mu \in \Bigl\{-\frac{\pi}{2}, \frac{\pi}{2}\Bigr\},\lambda \in \{1,2\},\sigma \in \{1,2\}$. We report the results for von Mises process in Table \ref{tab:vmp_estimator_005} for $\Delta t=0.005$, in Table \ref{tab:vmp_estimator_05} for  $\Delta t=0.05$ and Table \ref{tab:vmp_estimator_5} for $\Delta t=0.5$. We measure the bias and the standard deviation of the difference between the true parameter and estimated parameter for $\lambda$ and  $\sigma$. Since,  $\mu$ is a circular variable we measure its bias and concentration  by,

 \begin{align*}
	 \text{Bias}&=\text{atan2}(S,C)\\
	 \text{Concentration}&=\sqrt{S^2+C^2}
\end{align*}

where $S=\sum_{k=1}^{100}\sin(\mu-\hat{\mu}_k)$ and  $C=\sum_{k=1}^{100}\cos(\mu-\hat{\mu}_k)$.

A good estimator of $\mu$ should possess a low value of Bias, i.e. close to 0 and a high value of Concentration , i.e. close to 1.

We observe that, bias in $\hat{\lambda}$ and $\hat{\mu}$ is larger for smaller $n$, improving as $n$ increases. $\hat{\sigma}$ for the von Mises process also shows minimal bias. As $\Delta t$ grows larger and $n$ is small, the drift parameter $\hat{\lambda}$ shows moderate estimation error. This is because there are fewer data points to capture the mean-reversion behavior well. The bias of $\hat{\mu}$ shrinks with increasing $n$, and the concentration metric nears 1, indicating accurate estimates of $\mu$.

\begin{table}[H]
\scalebox{0.8}{%
\centering
\begin{tabular}{|c|c|c|c|c|c|c|c|c|c|} 
\hline
$\mu$   & $\lambda$ & $\sigma$ & $n$     & $\mathbb{E}[\lambda-\hat{\lambda}]$ & $\sqrt{\text{Var}[\lambda-\hat{\lambda}]}$ & $\mathbb{E}[\sigma-\hat{\sigma}]$ & $\sqrt{\text{Var}[\sigma-\hat{\sigma}]}$ & Bias of $\hat{\mu}$ & Concentration of $\mu-\hat{\mu}$ \\ 
\hline
-$\pi/2$ & 1      & 1     & 1000  & -1.385       & 1.501       & 0.002       & 0.023      & -0.179   & 0.787    \\ 
\hline
-$\pi/2$ & 1      & 2     & 1000  & -1.094       & 1.232       & 0.002       & 0.042      & -0.291   & 0.471    \\ 
\hline
-$\pi/2$ & 2      & 1     & 1000  & -1.346       & 1.562       & -0.001      & 0.023      & -0.087   & 0.947    \\ 
\hline
-$\pi/2$ & 2      & 2     & 1000  & -1.053       & 1.463       & 0.000       & 0.042      & 0.021    & 0.784    \\ 
\hline
$\pi/2$  & 1      & 1     & 1000  & -1.400       & 1.452       & 0.000       & 0.022      & 0.243    & 0.757    \\ 
\hline
$\pi/2$  & 1      & 2     & 1000  & -1.051       & 1.367       & 0.001       & 0.048      & 0.072    & 0.514    \\ 
\hline
$\pi/2$  & 2      & 1     & 1000  & -1.384       & 1.433       & 0.000       & 0.023      & 0.074    & 0.938    \\ 
\hline
$\pi/2$  & 2      & 2     & 1000  & -1.050       & 1.511       & -0.002      & 0.048      & 0.073    & 0.803    \\ 
\hline
-$\pi/2$ & 1      & 1     & 5000  & -0.692       & 0.533       & -0.001      & 0.009      & -0.064   & 0.954    \\ 
\hline
-$\pi/2$ & 1      & 2     & 5000  & -0.360       & 0.545       & -0.004      & 0.021      & -0.212   & 0.855    \\ 
\hline
-$\pi/2$ & 2      & 1     & 5000  & -1.014       & 0.714       & -0.002      & 0.009      & -0.022   & 0.991    \\ 
\hline
-$\pi/2$ & 2      & 2     & 5000  & -0.316       & 0.678       & -0.001      & 0.020      & -0.072   & 0.969    \\ 
\hline
$\pi/2$  & 1      & 1     & 5000  & -0.612       & 0.575       & -0.002      & 0.010      & 0.028    & 0.955    \\ 
\hline
$\pi/2$  & 1      & 2     & 5000  & -0.381       & 0.557       & -0.004      & 0.019      & 0.185    & 0.822    \\ 
\hline
$\pi/2$  & 2      & 1     & 5000  & -1.018       & 0.776       & -0.003      & 0.010      & 0.012    & 0.990    \\ 
\hline
$\pi/2$  & 2      & 2     & 5000  & -0.446       & 0.644       & -0.004      & 0.021      & 0.067    & 0.966    \\ 
\hline
-$\pi/2$ & 1      & 1     & 10000 & -0.578       & 0.355       & -0.001      & 0.008      & -0.003   & 0.983    \\ 
\hline
-$\pi/2$ & 1      & 2     & 10000 & -0.246       & 0.440       & 0.001       & 0.015      & -0.273   & 0.939    \\ 
\hline
-$\pi/2$ & 2      & 1     & 10000 & -0.999       & 0.510       & -0.001      & 0.007      & -0.002   & 0.996    \\ 
\hline
-$\pi/2$ & 2      & 2     & 10000 & -0.298       & 0.519       & -0.003      & 0.014      & -0.076   & 0.984    \\ 
\hline
$\pi/2$  & 1      & 1     & 10000 & -0.609       & 0.434       & -0.001      & 0.008      & 0.022    & 0.981    \\ 
\hline
$\pi/2$  & 1      & 2     & 10000 & -0.161       & 0.402       & -0.003      & 0.016      & 0.216    & 0.910    \\ 
\hline
$\pi/2$  & 2      & 1     & 10000 & -0.946       & 0.442       & -0.002      & 0.007      & 0.016    & 0.996    \\ 
\hline
$\pi/2$  & 2      & 2     & 10000 & -0.341       & 0.487       & -0.003      & 0.013      & 0.088    & 0.986    \\
\hline
\end{tabular}%

}
\caption{Simulation results for the von Mises process for $\Delta t=0.005$. Columns report the bias and the standard deviation of $\lambda-\hat{\lambda}$ and $\sigma-\hat{\sigma}$. We also report the bias and concentration for  $\hat{\mu}$, see section \ref{sec:sim_study} for their definitions.}
\label{tab:vmp_estimator_005}
\end{table}

\begin{table}[H]
\scalebox{0.8}{%
\centering
\begin{tabular}{|c|c|c|c|c|c|c|c|c|c|c|} 
\hline
$\mu$   & $\lambda$ & $\sigma$ & $n$     & $\mathbb{E}[\lambda-\hat{\lambda}]$ & $\sqrt{\text{Var}[\lambda-\hat{\lambda}]}$ & $\mathbb{E}[\sigma-\hat{\sigma}]$ & $\sqrt{\text{Var}[\sigma-\hat{\sigma}]}$ & Bias of $\hat{\mu}$ & Concentration of $\mu-\hat{\mu}$ \\ 
\hline
-$\pi/2$ & 1 & 1 & 1000  & -0.131 & 0.254 & -0.002 & 0.022 & -0.059 & 0.989  \\ 
\hline
-$\pi/2$ & 1 & 2 & 1000  & -0.226 & 0.420 & -0.004 & 0.044 & -0.294 & 0.926  \\ 
\hline
-$\pi/2$ & 2 & 1 & 1000  & -0.156 & 0.311 & -0.021 & 0.023 & -0.010 & 0.997  \\ 
\hline
-$\pi/2$ & 2 & 2 & 1000  & -0.308 & 0.497 & -0.022 & 0.042 & -0.113 & 0.986  \\ 
\hline
$\pi/2$  & 1 & 1 & 1000  & -0.154 & 0.247 & -0.007 & 0.022 & 0.032  & 0.989  \\ 
\hline
$\pi/2$  & 1 & 2 & 1000  & -0.341 & 0.460 & -0.008 & 0.047 & 0.198  & 0.927  \\ 
\hline
$\pi/2$  & 2 & 1 & 1000  & -0.145 & 0.325 & -0.020 & 0.022 & 0.006  & 0.997  \\ 
\hline
$\pi/2$  & 2 & 2 & 1000  & -0.204 & 0.487 & -0.021 & 0.046 & 0.095  & 0.984  \\ 
\hline
-$\pi/2$ & 1 & 1 & 5000  & -0.094 & 0.108 & -0.009 & 0.009 & -0.020 & 0.997  \\ 
\hline
-$\pi/2$ & 1 & 2 & 5000  & -0.159 & 0.223 & -0.007 & 0.021 & -0.253 & 0.988  \\ 
\hline
-$\pi/2$ & 2 & 1 & 5000  & -0.079 & 0.119 & -0.021 & 0.009 & 0.000  & 0.999  \\ 
\hline
-$\pi/2$ & 2 & 2 & 5000  & -0.276 & 0.235 & -0.023 & 0.021 & -0.105 & 0.997  \\ 
\hline
$\pi/2$  & 1 & 1 & 5000  & -0.098 & 0.104 & -0.008 & 0.011 & 0.019  & 0.997  \\ 
\hline
$\pi/2$  & 1 & 2 & 5000  & -0.174 & 0.188 & -0.006 & 0.019 & 0.252  & 0.988  \\ 
\hline
$\pi/2$  & 2 & 1 & 5000  & -0.084 & 0.136 & -0.023 & 0.009 & 0.006  & 0.999  \\ 
\hline
$\pi/2$  & 2 & 2 & 5000  & -0.230 & 0.202 & -0.022 & 0.018 & 0.110  & 0.997  \\ 
\hline
-$\pi/2$ & 1 & 1 & 10000 & -0.090 & 0.074 & -0.010 & 0.007 & -0.022 & 0.999  \\ 
\hline
-$\pi/2$ & 1 & 2 & 10000 & -0.153 & 0.148 & -0.007 & 0.013 & -0.248 & 0.995  \\ 
\hline
-$\pi/2$ & 2 & 1 & 10000 & -0.058 & 0.080 & -0.021 & 0.007 & -0.001 & 1.000  \\ 
\hline
-$\pi/2$ & 2 & 2 & 10000 & -0.228 & 0.139 & -0.023 & 0.013 & -0.107 & 0.999  \\ 
\hline
$\pi/2$  & 1 & 1 & 10000 & -0.087 & 0.076 & -0.009 & 0.008 & 0.028  & 0.999  \\ 
\hline
$\pi/2$  & 1 & 2 & 10000 & -0.140 & 0.161 & -0.005 & 0.013 & 0.270  & 0.995  \\ 
\hline
$\pi/2$  & 2 & 1 & 10000 & -0.073 & 0.089 & -0.023 & 0.007 & 0.004  & 1.000  \\ 
\hline
$\pi/2$  & 2 & 2 & 10000 & -0.239 & 0.144 & -0.023 & 0.015 & 0.105  & 0.999  \\
\hline
\end{tabular}%
}
\caption{Simulation results for the von Mises process for $\Delta t=0.05$. Columns report the bias and the standard deviation of $\lambda-\hat{\lambda}$ and $\sigma-\hat{\sigma}$. We also report the bias and concentration for  $\hat{\mu}$, see section \ref{sec:sim_study} for their definitions.}
\label{tab:vmp_estimator_05}
\end{table}

\begin{table}[H]
\scalebox{0.8}{%
\centering
\begin{tabular}{|c|c|c|c|c|c|c|c|c|c|} 
\hline
$\mu$   & $\lambda$ & $\sigma$ & $n$     & $\mathbb{E}[\lambda-\hat{\lambda}]$ & $\sqrt{\text{Var}[\lambda-\hat{\lambda}]}$ & $\mathbb{E}[\sigma-\hat{\sigma}]$ & $\sqrt{\text{Var}[\sigma-\hat{\sigma}]}$ & Bias of $\hat{\mu}$ & Concentration of $\mu-\hat{\mu}$ \\ 
\hline
-$\pi/2$ & 1      & 1     & 1000  & -0.079       & 0.099       & -0.093      & 0.028      & -0.056   & 0.999    \\ 
\hline
-$\pi/2$ & 1      & 2     & 1000  & 0.106        & 0.168       & 0.041       & 0.034      & -0.364   & 0.984    \\ 
\hline
-$\pi/2$ & 2      & 1     & 1000  & -0.273       & 0.108       & -0.286      & 0.032      & -0.019   & 1.000    \\ 
\hline
-$\pi/2$ & 2      & 2     & 1000  & 0.441        & 0.191       & -0.082      & 0.042      & -0.175   & 0.995    \\ 
\hline
$\pi/2$  & 1      & 1     & 1000  & -0.080       & 0.078       & -0.095      & 0.026      & 0.040    & 0.999    \\ 
\hline
$\pi/2$  & 1      & 2     & 1000  & 0.113        & 0.164       & 0.035       & 0.039      & 0.316    & 0.976    \\ 
\hline
$\pi/2$  & 2      & 1     & 1000  & -0.259       & 0.100       & -0.284      & 0.030      & 0.023    & 1.000    \\ 
\hline
$\pi/2$  & 2      & 2     & 1000  & 0.414        & 0.198       & -0.077      & 0.033      & 0.190    & 0.994    \\ 
\hline
-$\pi/2$ & 1      & 1     & 5000  & -0.083       & 0.037       & -0.096      & 0.011      & -0.054   & 1.000    \\ 
\hline
-$\pi/2$ & 1      & 2     & 5000  & 0.161        & 0.078       & 0.040       & 0.014      & -0.340   & 0.997    \\ 
\hline
-$\pi/2$ & 2      & 1     & 5000  & -0.269       & 0.044       & -0.287      & 0.014      & -0.016   & 1.000    \\ 
\hline
-$\pi/2$ & 2      & 2     & 5000  & 0.433        & 0.075       & -0.079      & 0.017      & -0.193   & 0.999    \\ 
\hline
$\pi/2$  & 1      & 1     & 5000  & -0.071       & 0.036       & -0.096      & 0.012      & 0.056    & 1.000    \\ 
\hline
$\pi/2$  & 1      & 2     & 5000  & 0.144        & 0.081       & 0.038       & 0.017      & 0.340    & 0.997    \\ 
\hline
$\pi/2$  & 2      & 1     & 5000  & -0.270       & 0.056       & -0.287      & 0.013      & 0.018    & 1.000    \\ 
\hline
$\pi/2$  & 2      & 2     & 5000  & 0.426        & 0.087       & -0.079      & 0.017      & 0.189    & 0.999    \\ 
\hline
-$\pi/2$ & 1      & 1     & 10000 & -0.074       & 0.026       & -0.095      & 0.008      & -0.048   & 1.000    \\ 
\hline
-$\pi/2$ & 1      & 2     & 10000 & 0.152        & 0.055       & 0.038       & 0.012      & -0.345   & 0.998    \\ 
\hline
-$\pi/2$ & 2      & 1     & 10000 & -0.265       & 0.037       & -0.286      & 0.010      & -0.018   & 1.000    \\ 
\hline
-$\pi/2$ & 2      & 2     & 10000 & 0.435        & 0.054       & -0.079      & 0.012      & -0.192   & 1.000    \\ 
\hline
$\pi/2$  & 1      & 1     & 10000 & -0.076       & 0.025       & -0.096      & 0.008      & 0.047    & 1.000    \\ 
\hline
$\pi/2$  & 1      & 2     & 10000 & 0.152        & 0.058       & 0.038       & 0.010      & 0.349    & 0.998    \\ 
\hline
$\pi/2$  & 2      & 1     & 10000 & -0.261       & 0.029       & -0.287      & 0.009      & 0.019    & 1.000    \\ 
\hline
$\pi/2$  & 2      & 2     & 10000 & 0.431        & 0.063       & -0.080      & 0.013      & 0.189    & 0.999    \\
\hline
\end{tabular}%
}
\caption{Simulation results for the von Mises process for $\Delta t=0.5$. Columns report the bias and the standard deviation of $\lambda-\hat{\lambda}$ and $\sigma-\hat{\sigma}$. We also report the bias and concentration for  $\hat{\mu}$, see section \ref{sec:sim_study} for their definitions.}
\label{tab:vmp_estimator_5}
\end{table}

Overall the simulation study confirms that the quadratic-variation-based estimator for $\sigma$ is accurate for both the circular Brownian motion and the von Mises process. The approximate MLE strategy for the von Mises process (using the derived transition density approximation) performs well, with parameter estimates converging as $n$ increases and as $\Delta t$ becomes smaller.

\subsection{Real data analysis}

We now illustrate an application of the estimation of the above models to wind direction data. The following data is obtained from a wind-power farm named ``Sotavento Galicia", located in Spain. We also refer to it as Sotavento Wind Farm in this paper. Sotavento Wind Farm makes available the data for wind direction at their farm, \cite{sotavento_historical_data}. The data from the Sotavento wind farm has been previously analyzed in the literature, see \cite{khazaei2022high}. The analysis of wind data is relevant for evaluating the economic and technical viability of wind power projects. 

\begin{figure}[H]
\includegraphics[scale=0.7]{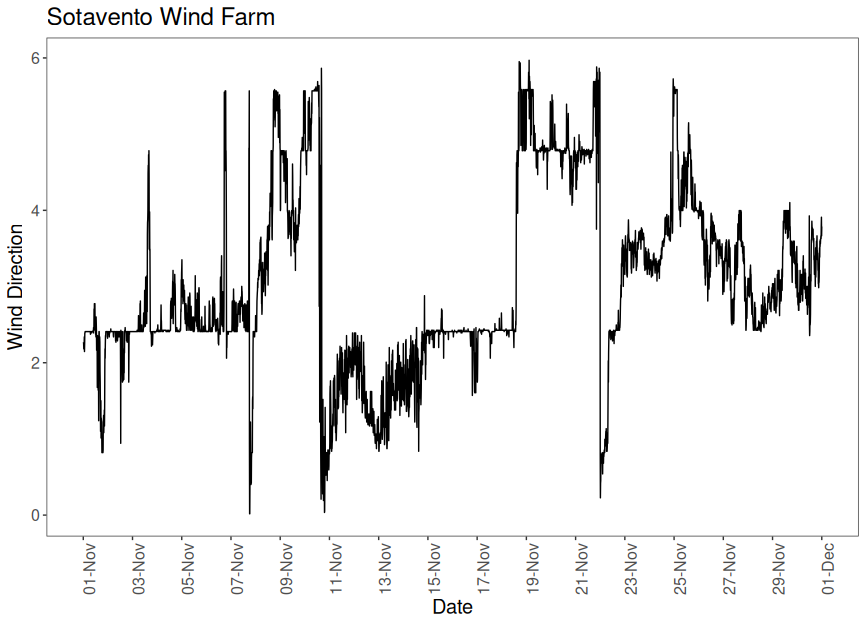}
\caption{Wind direction data in radians from 1 November 2024 to 30 November 2024 at a 10 minute frequency.}
\label{fig:sotavento}
\end{figure}

In our analysis we obtain wind direction data from Sotavento Wind farm at a 10 minute frequency from 1 November 2024 to 30 November 2024. We plot the time series of this wind direction data in Figure \ref{fig:sotavento}.

We fit both the circular Brownian motion and the von Mises process to this data. Let $\theta_t$ denote the wind-direction at time  $t$. The circular Brownian motion model is given by,
 $$
 d\theta_t=\sigma_{\text{bm}}dW_t
$$

and the von Mises process model is given by,

$$
d\theta_t=-\lambda_{\text{vm}}\sin(\theta_t-\mu_{\text{vm}})dt +\sigma_{\text{vm}}dW_t
$$

As per the notation above $t_{i+1}-t_i=10\text{ minutes}, \forall i$. Therefore  $\Delta t=t_{i+1}-t_i=\frac{10}{1440}$ days.

We report that $\hat{\sigma}_{\text{bm}}=2.96$. For the von Mises process, the estimates are $\hat{\sigma}_{\text{vm}}=2.96, \hat{\lambda}_{\text{vm}}=4.60$ and  $\hat{\mu}_{\text{vm}}=2.61$. The estimated  $\hat{\sigma}_{\text{bm}}$ and  $\hat{\sigma}_{\text{vm}}$ are identical because both were obtained via the quadratic variation estimator. We also obtain the 95\% confidence region for the parameters using parametric bootstrap. We describe the parametric bootstrap procedure for the stochastic correlation model in Algorithm \ref{bootstrap}.

The 95\% confidence region for $\hat{\sigma}_{\text{bm}}$ is $(2.90,3.02)$ and for  $\hat{\sigma}_{\text{vm}}$ is $(2.89,3.01)$. We note that under both models the estimated volatility coefficient and their intervals are nearly identical. The 95\% confidence region for $\hat{\mu}_{\text{vm}}$ is  $(2.11,3.11)$. The 95\% confidence interval for  $\hat{\lambda}_{\text{vm}}$ is  $(3.25,6.74)$. This indicates a mean-reverting behavior of the wind-direction at the Sotavento farm.

\section{Stochastic correlation models in finance}\label{sec: stochcorr}

In this section we propose to model stochastic correlation in financial markets using the circular diffusions described above. We also consider the estimation of such a model.

Let $S^{(1)}_t$ and $S^{(2)}_t$ describe the price of two assets at time $t$. We assume that the price dynamics of $S^{(1)}_t$ and  $S^{(2)}_t$  follow a correlated geometric Brownian motion. Geometric Brownian motion is a benchmark stock price model due to its analytical tractability, other suitable price processes can also be considered. We consider two cases- one where the volatility is constant and the other where the volatility is stochastic, with correlation being stochastic in both cases. 

\subsection{Stochastic correlation with constant volatility}

We describe the dynamics below,

\begin{align}\label{stochcorr_model}
	dS^{(1)}_t&=\mu^{(1)}S^{(1)}_t dt +\sigma^{(1)}S^{(1)}_t dB^{(1)}_t\nonumber\\
	dS^{(2)}_t&=\mu^{(2)}S^{(2)}_t dt +\sigma^{(2)}S^{(2)}_t\left(\rho_t dB^{(1)}_t +\sqrt{1-\rho_t^2}dB^{(2)}_t\right)
\end{align}

where $\rho_t=\cos\theta_t$, with two possible choices for modelling $\theta_t$ being the circular Brownian motion ($d\theta_t=\sigma dB^{(3)}_t$) and the von Mises process ($d\theta_t=-\lambda\sin(\theta_t-\mu)dt+\sigma dB^{(3)}_t$). Here $B^{(1)}_t,B^{(2)}_t,B^{(3)}_t$ are mutually independent standard Brownian motion process. In this specification, note that the correlation process is independent from the price process. This assumption can be relaxed at the cost of some loss in tractability during estimation. 
\begin{theorem}\label{thm:strong}
    The strong solution to SDEs (\ref{stochcorr_model}) model exists.
\end{theorem}
\begin{proof}
	See the appendix for a proof.
\end{proof}

\subsection{Stochastic correlation with stochastic volatility}

To account for the time dependency in volatility in geometric Brownian motion based model for stock prices, \cite{heston1993closed} introduced a stochastic volatility model where the squared volatility of the asset price follows a Cox-Ingersoll-Ross process (\cite{Cox1985}). To incorporate the effect of stochastic correlation we consider a bivariate version of the Heston model,
\begin{align*}
	dS_t^{(1)}&=\mu^{(1)}S_t^{(1)}dt+\sqrt{v_t^{(1)}}S_t^{(1)}dB_t^{S^{(1)}}\\
	dv_t^{(1)}&=-\kappa^{(1)}(v_t^{(1)}-\theta^{(1)})dt+\sigma^{(1)}\sqrt{v_t^{(1)}}dB_t^{v^{(1)}}\\
	dS_t^{(2)}&=\mu^{(2)}S_t^{(2)}dt+\sqrt{v_t^{(2)}}S_t^{(2)}dB_t^{S^{(2)}}\\
	dv_t^{(2)}&=-\kappa^{(2)}(v_t^{(2)}-\theta^{(2)})dt+\sigma^{(2)}\sqrt{v_t^{(2)}}dB_t^{v^{(2)}}
\end{align*}

In the Heston model, the asset price and its volatility are correlated. This is given by, $dB_t^{S^{(1)}}dB_t^{v^{(1)}}=\rho^{(1)}dt$ and $dB_t^{S^{(2)}}dB_t^{v^{(2)}}=\rho^{(2)}dt$, where $\rho^{(1)}$ is referred to as leverage of  $S_t^{(1)}$ and  $\rho^{(2)}$ as the leverage of  $S_t^{(2)}$. It's further assumed that both the asset prices are also stochastically correlated, i.e., $dB_t^{S^{(1)}}dB_t^{S^{(2)}}=\rho_t dt$.

We consider an alternate decomposition of the driving Brownian motion in terms of independent Brownian motions as follows,
\begin{align*}
	dB_t^{S^{(1)}}&=dB_t^{(1)}\\
	dB_t^{v^{(1)}}&=\rho^{(1)}dB_t^{(1)}+\sqrt{1-(\rho^{(1)})^2}dB_t^{(2)}\\
	dB_t^{S^{(2)}}&=\rho_t dB_t^{(1)}+\sqrt{1-\rho_t^2}dB_t^{(3)}\\
	dB_t^{v^{(2)}}&=\rho^{(2)}\rho_t dB_t^{(1)}+\rho^{(2)}\sqrt{1-\rho_t^2}dB_t^{(3)}+\sqrt{1-(\rho^{(2)})^2}dB_t^{(4)}
\end{align*}

here $B_t^{(i)},i=1,2,3,4$ are mutually independent Brownian motions. This decomposition allows us to approximately simulate from the bivariate Heston model through its Euler-Maruyama discretization.

\subsection{Estimation for the constant volatility specification}\label{sec:constvol_est}

In the stochastic correlation model the instantaneous correlation
$\rho_t$ is latent. Estimation of models with latent states is often
approached via filtering methods; see, for example, \cite{lautier2003filtering}.
Here we exploit the fact that, conditional on a realization of the correlation
path, the log-returns are Gaussian and admit a closed-form likelihood.
This allows us to estimate both the model parameters and a most likely
realization of the correlation path by joint maximization.

Assume we observe prices on an equidistant grid $t_i=i\Delta t$,
$i=0,1,\dots,n$, with $T=n\Delta t$. Denote the observations by
$$
\mathbf S^{(1)}=\{S^{(1)}_{t_0},S^{(1)}_{t_1},\ldots,S^{(1)}_{t_n}\},
\qquad
\mathbf S^{(2)}=\{S^{(2)}_{t_0},S^{(2)}_{t_1},\ldots,S^{(2)}_{t_n}\}.
$$

The (unobserved) correlation path is written as
$$
\boldsymbol\rho=\{\rho_{t_0},\rho_{t_1},\ldots,\rho_{t_n}\},
\qquad
\rho_{t_i}=\cos(\theta_{t_i}),
$$
where $\theta_t$ follows either the circular Brownian motion or the von Mises
diffusion introduced earlier. Throughout, we assume that the correlation
process is independent of the Brownian motions driving the asset prices.

Let $X^{(j)}_{t_i}=\log S^{(j)}_{t_i}$ denote log-prices. Under the constant
volatility specification, conditional on $\rho_{t_i}$ the vector of
log-increments
$$
\Delta X_i=
\begin{pmatrix}
X^{(1)}_{t_{i+1}}-X^{(1)}_{t_i}\\[2pt]
X^{(2)}_{t_{i+1}}-X^{(2)}_{t_i}
\end{pmatrix}
$$

is Gaussian with mean vector
$$
m=
\begin{pmatrix}
(\mu^{(1)}-\tfrac12(\sigma^{(1)})^2)\Delta t\\[2pt]
(\mu^{(2)}-\tfrac12(\sigma^{(2)})^2)\Delta t
\end{pmatrix}
$$

and covariance matrix
$$
\Sigma_i
=
\Delta t\,
\begin{pmatrix}
(\sigma^{(1)})^2 & \sigma^{(1)}\sigma^{(2)}\rho_{t_i}\\[2pt]
\sigma^{(1)}\sigma^{(2)}\rho_{t_i} & (\sigma^{(2)})^2
\end{pmatrix}.
$$

Therefore the conditional log-likelihood of the asset prices given
$\boldsymbol\rho$ can be written as
\begin{equation}\label{eq:cond_ll_constvol}
\tilde L(\mu^{(1)},\sigma^{(1)},\mu^{(2)},\sigma^{(2)};\mathbf S^{(1)},\mathbf S^{(2)},\boldsymbol\rho)
=
-\frac12\sum_{i=0}^{n-1}
\left\{
(\Delta X_i-m)^\top\Sigma_i^{-1}(\Delta X_i-m)
+\log|\Sigma_i|
\right\},
\end{equation}
up to an additive constant independent of the parameters.

Let $p(\cdot\mid\cdot)$ denote the transition density of the chosen angular
process for $\theta_t$ over time step $\Delta t$. The correlation is defined as
$\rho_t=\cos(\theta_t)$, and the mapping $\rho=\cos\theta$ from $[0,2\pi)$ to
$(-1,1)$ is two-to-one. For any $\rho\in(-1,1)$, the two preimages are
$$
\theta^{(1)}=\arccos(\rho),
\qquad
\theta^{(2)}=2\pi-\arccos(\rho).
$$

In general, the induced transition density of $\rho_{t_{i+1}}$ given
$\rho_{t_i}$ is obtained by summing the contributions from both branches of the
inverse transformation. In practice, we approximate this sum by evaluating the
angular transition density at the principal branch and doubling the
contribution. This approximation is motivated by the observation that the diffusion transition densities are locally normal for sufficiently small $\Delta t$, the
dominant contribution arises from the nearest preimage
$\arccos(\rho_{t_{i+1}})$, while the second branch contributes an equal amount
to first order. Under this approximation, the induced log-likelihood
contribution is taken to be
$$
\log p\!\left(\arccos(\rho_{t_{i+1}})\mid \arccos(\rho_{t_i})\right)
+
\log\!\left(\frac{2}{\sqrt{1-\rho_{t_{i+1}}^2}}\right).
$$

This expression is used in all subsequent likelihood evaluations and numerical
optimization.

Since the correlation process is independent of the asset-price innovations,
the joint log-likelihood for parameters and latent correlation path is
\begin{align}
L(\vartheta;\mathbf S^{(1)},\mathbf S^{(2)},\boldsymbol\rho)
&=
\tilde L(\mu^{(1)},\sigma^{(1)},\mu^{(2)},\sigma^{(2)};\mathbf S^{(1)},\mathbf S^{(2)},\boldsymbol\rho)
\nonumber\\
&\quad+
\sum_{i=0}^{n-1}
\left[
\log p\!\left(\arccos(\rho_{t_{i+1}})\mid \arccos(\rho_{t_i})\right)
+\log\!\left(\frac{2}{\sqrt{1-\rho_{t_{i+1}}^2}}\right)
\right],
\label{eq:joint_ll_constvol}
\end{align}
where
$$
\vartheta=
\bigl(\mu^{(1)},\sigma^{(1)},\mu^{(2)},\sigma^{(2)},\text{(correlation parameters)}\bigr).
$$

In particular, if $\theta_t$ is the circular Brownian motion, $p$ is given by
\eqref{cbm_tpd}, while for the von Mises diffusion we use the approximation in
\eqref{tpd_vm}.

Direct maximization of \eqref{eq:joint_ll_constvol} over the full latent path
$\boldsymbol\rho$ can lead to excessive inter-period variability in estimated
correlations. To stabilize the optimization and to enforce economically
plausible smoothness, we maximize a penalized criterion of the form
\begin{equation}\label{eq:pen_ll_constvol}
L^*(\vartheta,\boldsymbol\rho)
=
L(\vartheta;\mathbf S^{(1)},\mathbf S^{(2)},\boldsymbol\rho)
-
\lambda_1\,n\sum_{i=0}^{n-1}\bigl(\rho_{t_{i+1}}-\rho_{t_i}\bigr)^2
-
\lambda_2\,\kappa,
\end{equation}
where $\lambda_1\ge0$ controls smoothness of the correlation path and, in the
von Mises specification, $\kappa=2\lambda/\sigma^2$ is the concentration
parameter of the stationary von Mises distribution. The additional penalty in
$\kappa$ is motivated by numerical stability: evaluation of modified Bessel
functions in the von Mises transition density becomes unreliable for very large
$\kappa$. When the correlation process is the circular Brownian motion we set
$\lambda_2=0$.

We maximize $L^*$ numerically using the derivative-free constrained optimizer
COBYLA \cite{powell1998direct}, as implemented in the \texttt{NLopt} package
\cite{NLopt}. Since the correlation process is specified independently of the
asset-price dynamics, at each iteration we update the diffusion parameter of
the correlation process using a quadratic-variation-type estimator based on the
current iterate of $\boldsymbol\rho$; see \eqref{quad_var_est}. The resulting
procedure yields simultaneous estimates of the model parameters and a most
likely realization of the latent correlation path.

\subsection{Estimation for stochastic volatility specification}\label{sec:sv_stochcorr_est}

We now extend the constant-volatility specification of
Section \ref{sec:constvol_est} by allowing the marginal volatilities of both
assets to vary over time. The correlation path $\boldsymbol\rho$ is treated in
exactly the same way as before: it is latent, parametrized by
$\rho_{t_i}=\cos(\theta_{t_i})$ on the principal branch, and its transition
density is induced by the von Mises diffusion for $\theta_t$ together with the
Jacobian term. Consequently, the additional modelling and likelihood
contributions in this section arise solely from the introduction of latent
volatility paths.

In addition to $\boldsymbol\rho=\{\rho_{t_0},\ldots,\rho_{t_n}\}$, we introduce
two latent instantaneous volatility paths
$$
\mathbf v^{(1)}=\{v^{(1)}_{t_0},\ldots,v^{(1)}_{t_n}\},
\qquad
\mathbf v^{(2)}=\{v^{(2)}_{t_0},\ldots,v^{(2)}_{t_n}\},
$$

with $v^{(j)}_{t_i}>0$. (As in our implementation, $v^{(j)}$ denotes
instantaneous volatility rather than variance.)

Conditional on $(\rho_{t_{i+1}},v^{(1)}_{t_{i+1}},v^{(2)}_{t_{i+1}})$, the
log-return increment is approximated by a bivariate Gaussian density, with the
same structure as in \eqref{eq:cond_ll_constvol} but with time-varying
volatilities:
$$
\Delta X_i
=
\begin{pmatrix}
X^{(1)}_{t_{i+1}}-X^{(1)}_{t_i}\\[2pt]
X^{(2)}_{t_{i+1}}-X^{(2)}_{t_i}
\end{pmatrix}
\approx  N(m_i,\Sigma_i),
$$

where
$$
m_i=
\begin{pmatrix}
(\mu^{(1)}-\tfrac12 (v^{(1)}_{t_{i+1}})^2)\Delta t\\[2pt]
(\mu^{(2)}-\tfrac12 (v^{(2)}_{t_{i+1}})^2)\Delta t
\end{pmatrix},
\qquad
\Sigma_i
=
\Delta t\,
\begin{pmatrix}
(v^{(1)}_{t_{i+1}})^2 &
v^{(1)}_{t_{i+1}}v^{(2)}_{t_{i+1}}\rho_{t_{i+1}}\\[3pt]
v^{(1)}_{t_{i+1}}v^{(2)}_{t_{i+1}}\rho_{t_{i+1}}&
(v^{(2)}_{t_{i+1}})^2
\end{pmatrix}.
$$

Let $L_i^{S}$ denote the corresponding Gaussian log-density contribution.

We model each volatility path by
Cox-Ingersoll-Ross diffusion, and in the discretized likelihood we use the exact
one-step transition density for $v^{(j)}_{t_{i+1}}\mid v^{(j)}_{t_i}$. Denoting
this density by $p_{\mathrm{CIR}}(\cdot\mid\cdot)$, we write the corresponding
log-density contribution as
$$
L_i^{v^{(j)}}=\log p_{\mathrm{CIR}}\!\bigl(v^{(j)}_{t_{i+1}}\mid v^{(j)}_{t_i}\bigr),
\qquad j\in\{1,2\}.
$$

In our implementation the diffusion coefficients of the volatility processes are
updated from the current volatility paths using quadratic-variation-type
estimators,
$$
\hat\sigma^{(j)}_{\mathrm{vol}}
=
\left(
\frac{1}{n\Delta t}\sum_{i=0}^{n-1}\bigl(v^{(j)}_{t_{i+1}}-v^{(j)}_{t_i}\bigr)^2
\right)^{1/2},
\qquad j\in\{1,2\},
$$

analogously to the update of $\sigma_{\mathrm{vm}}$ from $\boldsymbol\rho$
described in Section~\ref{sec:constvol_est}.

The joint approximate log-likelihood is obtained by adding the volatility
transition contributions to the joint criterion of the constant-volatility
case,
$$
\mathcal L(\vartheta;\mathbf S^{(1)},\mathbf S^{(2)},\boldsymbol\rho,\mathbf v^{(1)},\mathbf v^{(2)})
=
\sum_{i=0}^{n-1}
\Bigl[
L_i^{S}
+
L_i^{v^{(1)}}
+
L_i^{v^{(2)}}
+
L_i^{\rho}
\Bigr],
$$

where $L_i^{\rho}$ denotes the correlation transition contribution (von Mises
transition density in the angular variable plus the Jacobian term, as in
Section~\ref{sec:constvol_est}). The finite-dimensional parameter vector is
$$
\vartheta=
\bigl(
\mu^{(1)},\mu^{(2)},
\lambda_{\mathrm{vm}},\mu_{\mathrm{vm}},
\kappa^{(1)},\theta^{(1)},
\kappa^{(2)},\theta^{(2)}
\bigr),
$$
with $\sigma_{\mathrm{vm}}$ and $\sigma^{(j)}_{\mathrm{vol}}$ updated from the
latent paths.

We employ the same smoothness penalty on $\boldsymbol\rho$ and the same
$\kappa_{\mathrm{vm}}$-penalty as in \eqref{eq:pen_ll_constvol}. Thus we maximize
a penalized criterion of the form
$$
\mathcal L^*
=
\mathcal L
-
\lambda_1\,n\sum_{i=0}^{n-1}\bigl(\rho_{t_{i+1}}-\rho_{t_i}\bigr)^2
-
\lambda_2\,\kappa_{\mathrm{vm}},
\qquad
\kappa_{\mathrm{vm}}=\frac{2\lambda_{\mathrm{vm}}}{\sigma_{\mathrm{vm}}^2}.
$$

Optimization is performed by the same class of derivative-free constrained
methods as in Section~\ref{sec:constvol_est}. In addition to the correlation
path, we now optimize over the volatility paths
$(\mathbf v^{(1)},\mathbf v^{(2)})$ subject to positivity constraints, together
with the finite-dimensional parameters in $\vartheta$. Initialization of
$\mathbf v^{(1)}$ and $\mathbf v^{(2)}$ is obtained from rolling-window
standard deviations of log-prices, while initialization of $\boldsymbol\rho$ is
as in the constant-volatility case.

The procedure returns the estimated latent paths
$\hat{\boldsymbol\rho}$, $\hat{\mathbf v}^{(1)}$, $\hat{\mathbf v}^{(2)}$, and
the associated parameter estimates.

\subsection{Bootstrap}

Since we work with penalized likelihood to estimate the model, the asymptotic normality of MLE is unavailable to quantify the uncertainty of the estimates. We perform parametric bootstrap (see p. 306, \cite{efron1994introduction}) to get confidence interval for the estimates of $\boldsymbol{\rho}$.

We recall a result from \cite{ducharme1985bootstrap,fisher1989bootstrap} regarding the choice of pivot for bootstrap for circular random variables. Let $\theta_1,\theta_2,\ldots,\theta_n$ denote the random sample of a circular random variable. Let $\hat{\alpha}$ be a circular estimate of a parameter from the sample. Define,  $U(\hat{\alpha},\alpha)=1-\cos(\hat{\alpha}-\alpha)$. Upon Bootstrap resampling, we obtain values  $U(\hat{\alpha},\alpha_1),U(\hat{\alpha},\alpha_2),\ldots,U(\hat{\alpha},\alpha_n)$. Then the $100\times q$-th percentile value of  $U$ is chosen corresponding to which we obtain  $\alpha_q$. Then $\theta_q=\hat{\alpha}-\alpha_q$, and therefore the confidence region for  $\hat{\alpha}$ is  $\hat{\alpha}\pm \alpha_q$.

Following this we propose $U^*(\hat{\rho},\rho)=1-\cos(\cos^{-1}(\hat{\rho})-\cos^{-1}(\rho))$ as the pivot for bootstrapping estimated correlation. Let  $U^*_q$ denote the  $100\times q$-th percentile value, then the $q-th$ confidence region for  $\hat{\rho}$ is $\hat{\rho}\pm U^*_q$. 

The bootstrap procedure is described in Algorithm \ref{bootstrap}.

\begin{algorithm}[H]
	\caption{Bootstrap confidence region for $\boldsymbol{\rho}$}
\label{bootstrap}
\begin{algorithmic}[1]
	\Statex \textbullet~\textbf{Input: } $\hat{\mu}^{(1)},\hat{\sigma}^{(1)},\hat{\mu}^{(2)},\hat{\sigma}^{(2)}, \boldsymbol{\hat{\rho}}=\{\hat{\rho}_0,\hat{\rho}_{\Delta t},\hat{\rho}_{2\Delta t},\ldots,\hat{\rho}_{T}\}$. Parameters of the circular diffusion process: circular Brownian motion ($\hat{\sigma}$) or von Mises process ($\hat{\lambda},\hat{\mu},\hat{\sigma}$).
\Statex \textbullet~\textbf{Parameters: }$N$: the number of bootstrap samples.
\Statex \textbullet~\textbf{Output: } Bootstrap confidence region for $\boldsymbol{\rho}$. 
\State Generate $N$ random samples from  $S^{(1)}_t$ and  $S^{(2)}_t$ at times  $t=\{0,\Delta t, 2\Delta t, \ldots, T\}$ using the estimated parameters as input. Then for each $i\in \{1,2,\ldots, N\}, \boldsymbol{S}^{(1)i}=\{S^{(1)i}_0,S^{(1)i}_{\Delta t},\ldots,S^{(1)i}_T\}$ and $\boldsymbol{S}^{(2)i}=\{S^{(2)i}_0,S^{(2)i}_{\Delta t},\ldots,S^{(2)i}_T\}$
\State For each $i$, we estimate the model for  $\boldsymbol{S}^{(1)i}$ and  $\boldsymbol{S}^{(2)i}$, corresponding to which we obtain an estimate  $\hat{\boldsymbol{\rho}}^i=\{\hat{\rho}_0^i,\hat{\rho}_{\Delta t}^i,\ldots, \hat{\rho}_T^i\}$. 
\State Then for each $t^*\in \{0,\Delta t,\ldots, T\}$, compute $U(\hat{\rho}_{t^*},\hat{\rho}^1_{t^*}),U(\hat{\rho}_{t^*},\hat{\rho}^2_{t^*}),\ldots,U(\hat{\rho}_{t^*},\hat{\rho}^N_{t^*})$. Let $U^{t^*}_q$ denote the  $100\times q$-th percentile value of these  $U$'s.
\State Then the $q$-th confidence region for  $\hat{\rho}_{t^*}$ is $\hat{\rho}_{t^*}\pm U^{t^*}_q$. 
\end{algorithmic}
\end{algorithm}

\begin{prop}
We note that the solution of the stochastic correlation model is given by, 
\begin{align*}
	S^{(1)}_t&=S^{(1)}_0\exp\left(\left(\mu^{(1)}-\frac{(\sigma^{(1)})^2}{2}\right)t+\sigma^{(1)}B_t^{(1)}\right)\\
	S^{(2)}_t&=S^{(2)}_0\exp\left(\left(\mu^{(2)}-\frac{(\sigma^{(2)})^2}{2}\right)t+\sigma^{(2)}(\rho_tB^{(1)}_t+\sqrt{1-\rho_t^2}B^{(2)}_t)\right)
\end{align*}
\end{prop}
We can plugin the estimates $\hat{\mu}^{(1)},\hat{\sigma}^{(1)},\hat{\mu}^{(2)},\hat{\sigma}^{(2)}, \boldsymbol{\hat{\rho}}$ to simulate from the model as is needed for Step 1 of the Bootstrap algorithm.

\section{Application to FX markets}\label{sec: application}

In this section we consider the application of the stochastic correlation model
to study the correlation dynamics of the foreign exchange (FX) and the stock
markets. Any cross-border trade carries a currency risk, i.e. since the
exchange rate is assumed to vary randomly in time, by the time terms of the
trade are fulfilled the currency rate could have changed enough to impact the
parties of the trade. To hedge such risks, financial markets participants enter
into specialized financial contracts like the quantity adjusting option, see
\cite{wystup2010}, wherein they can choose to lock-in the currency rate for some
pre-specified terms. As we noted in Section \ref{sec: intro} that there is work
towards understanding the pricing of the quantity adjusting options where the
currency and the equity are both correlated. Motivated by the interest in the
problem we seek to understand the variation in the correlation dynamics of the
currency market and the equity market around the world during COVID-19. It has
been noted that the shock of the COVID-19 pandemic was felt unlike other
previous shocks in the FX markets, \cite{bazan2021fx}.

We study the correlation dynamics between the following pairs of exchange rates
and equity indices- USD/EUR and S\&P500, GBP/USD and FTSE100, INR/USD and
NIFTY50, JPY/USD and NIKKEI225. Towards this we collect end-of-day closing data
for all the pairs for the year 2020 from Yahoo Finance. For each pair we fit
the stochastic correlation model under constant volatility and the
stochastic volatility extension described in Section \ref{sec:constvol_est} and \ref{sec:sv_stochcorr_est} respectively..
In both cases we consider two specifications for the correlation dynamics-the von Mises diffusion and the circular Brownian motion.

\begin{figure}[H]
\centering
\includegraphics[width=0.92\textwidth]{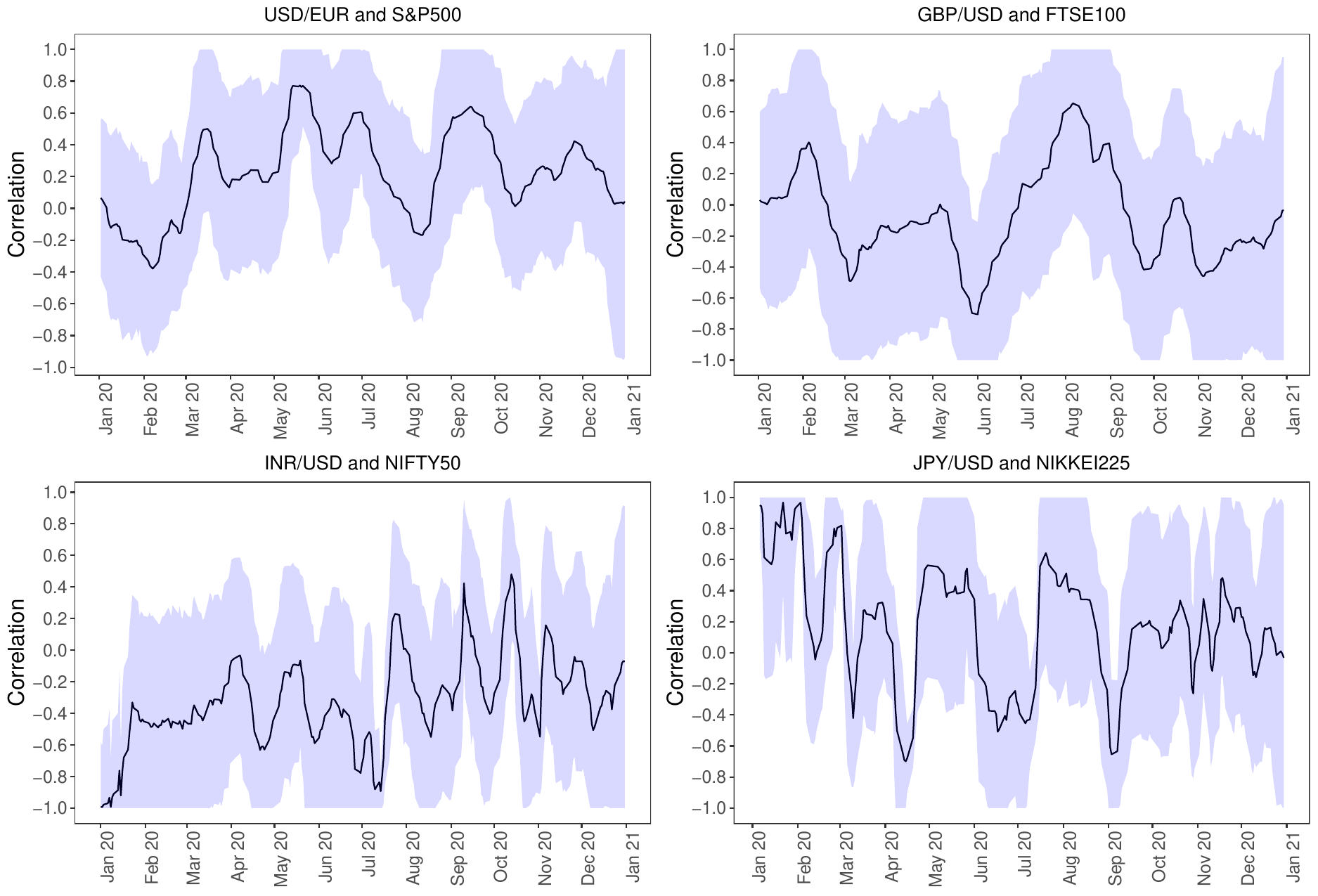}
\caption{Estimated correlation dynamics between USD/EUR and S\&P500, GBP/USD and
FTSE100, INR/USD and NIFTY50, and JPY/USD and NIKKEI225 for the year 2020 using
the von Mises process as the correlation model under the constant
volatility specification.}
\label{fig:vmp_constvol}
\end{figure}

\begin{figure}[H]
\centering
\includegraphics[width=0.92\textwidth]{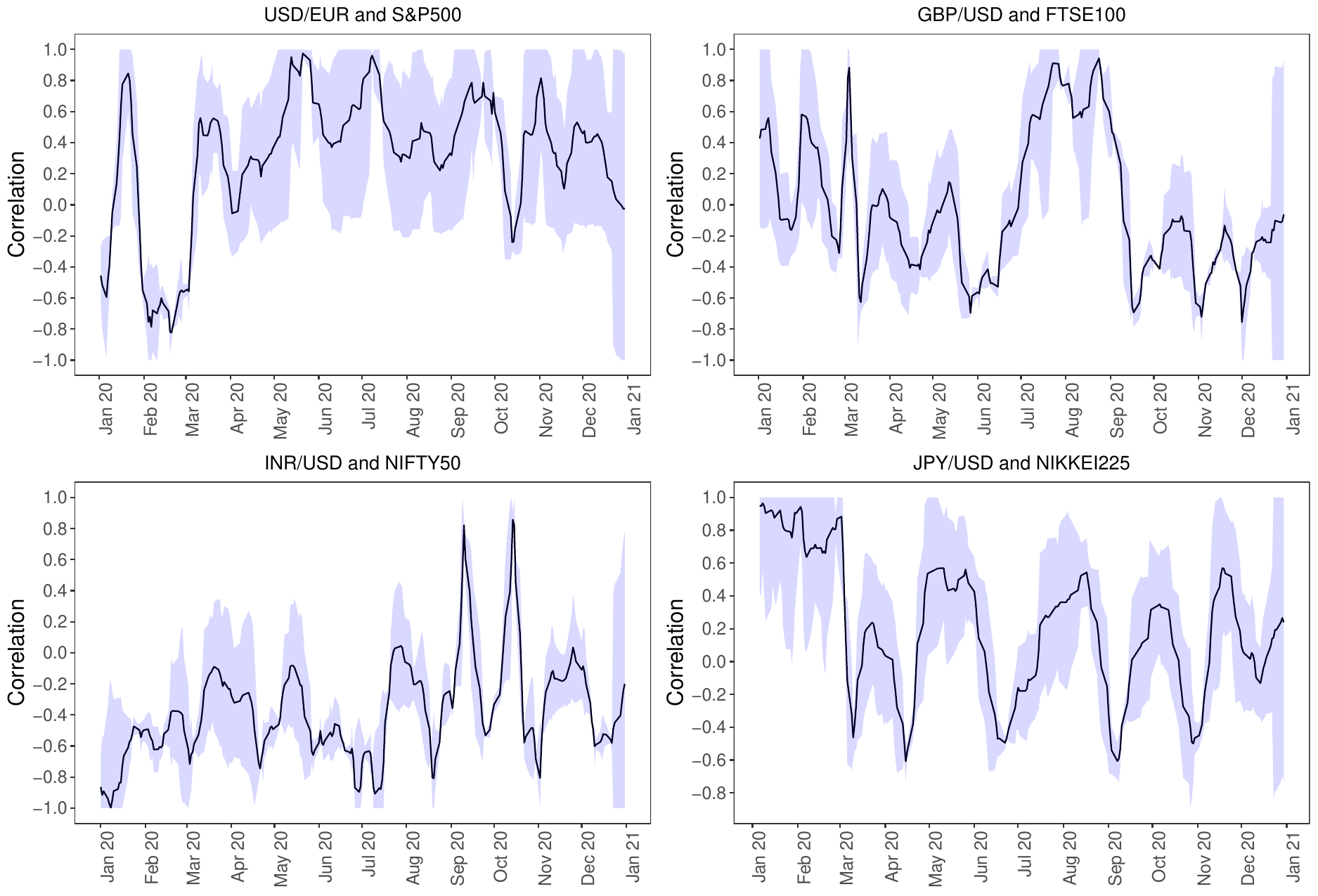}
\caption{Estimated correlation dynamics for the same four FX-equity pairs in
2020 using the von Mises process as the correlation model under the
stochastic volatility specification.}
\label{fig:vmp_stochvol}
\end{figure}

Figure \ref{fig:vmp_constvol} shows the estimated correlation dynamics for the
four FX--equity pairs when the correlation process is modelled by the von Mises
diffusion under the constant volatility specification. Figure~
\ref{fig:vmp_stochvol} shows the corresponding estimates under the
stochastic volatility specification.

Across all four pairs, the estimated correlation exhibits pronounced shocks
aligned with the onset of COVID-19 and subsequent restrictions. In particular,
for USD/EUR and S\&P500, an early-year spike followed by a rapid decline is
visible around January, coinciding with growing awareness of COVID-19 in the
USA. For the UK and India pairs, correlation spikes occur during March and
April, corresponding to the introduction of major restrictions. In most cases,
these shocks appear transitory rather than persistent. Allowing for stochastic
volatility preserves these broad patterns, while providing additional structure
during periods of heightened market turbulence.

\begin{figure}[H]
\centering
\includegraphics[width=0.92\textwidth]{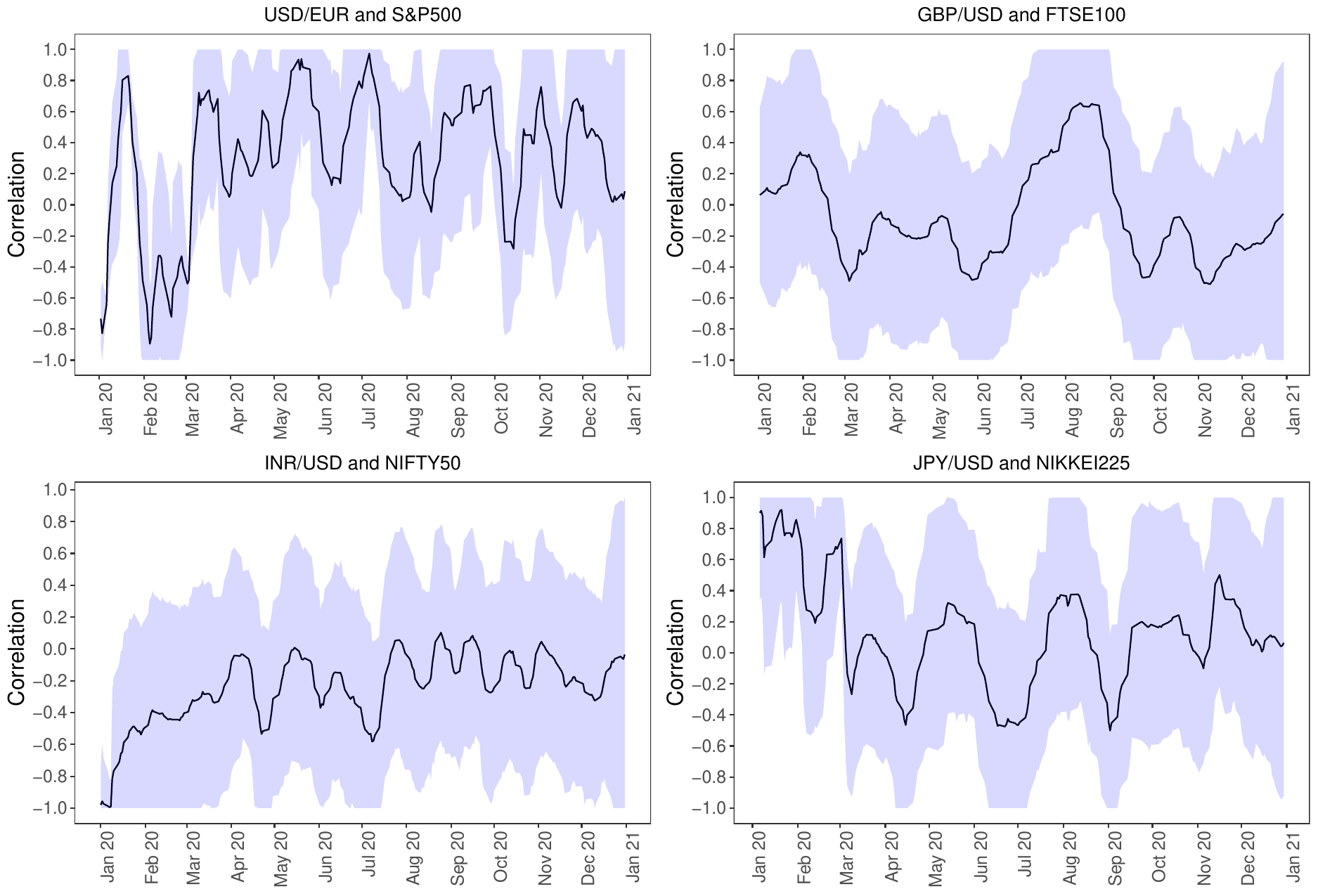}
\caption{Estimated correlation dynamics between USD/EUR and S\&P500, GBP/USD and
FTSE100, INR/USD and NIFTY50, and JPY/USD and NIKKEI225 for the year 2020 using
circular Brownian motion as the correlation model under the constant
volatility specification.}
\label{fig:cbm_constvol}
\end{figure}

\begin{figure}[H]
\centering
\includegraphics[width=0.92\textwidth]{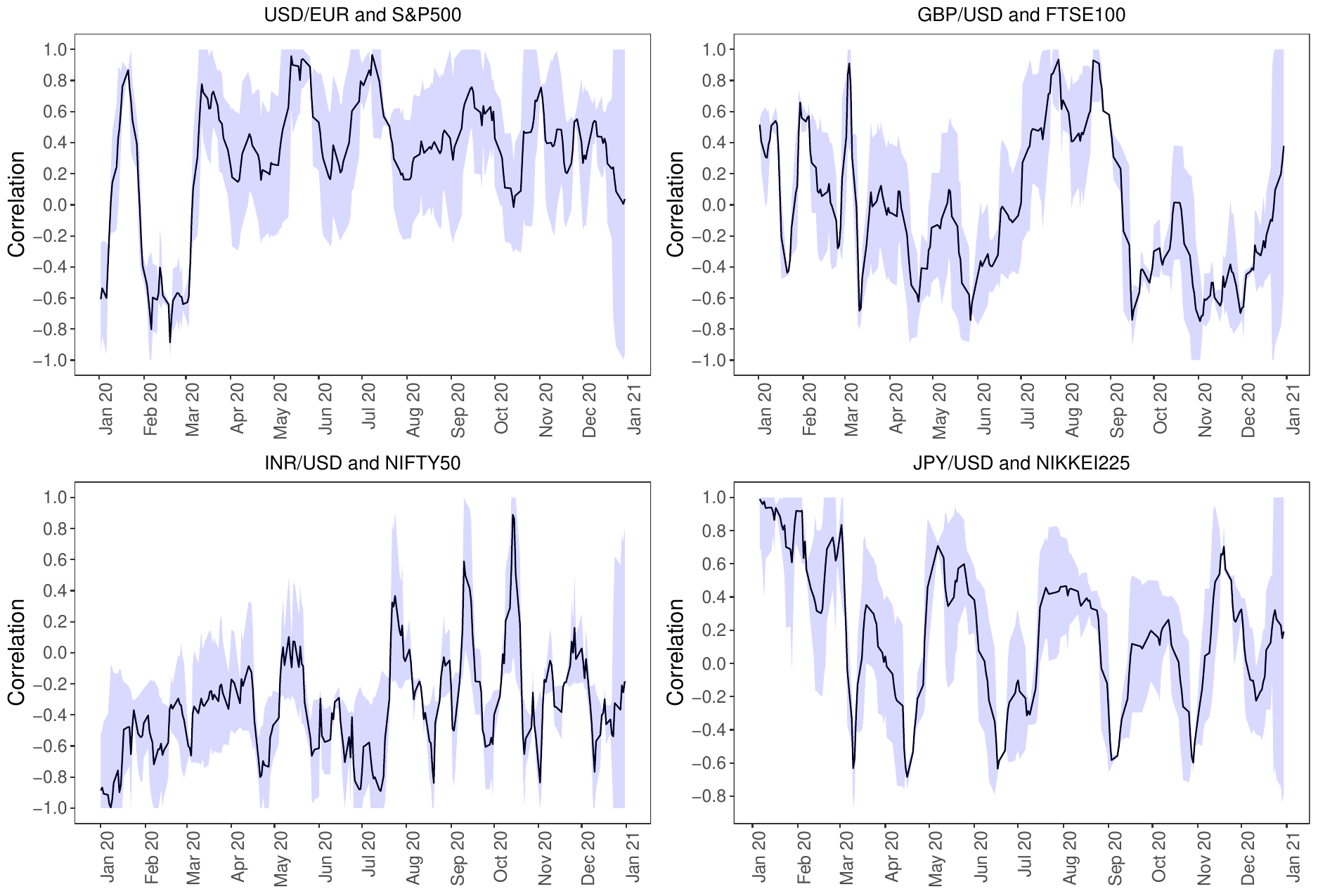}
\caption{Estimated correlation dynamics for the same four FX--equity pairs in
2020 using circular Brownian motion as the correlation model under the
stochastic volatility specification.}
\label{fig:cbm_stochvol}
\end{figure}

Figures \ref{fig:cbm_constvol} and \ref{fig:cbm_stochvol} report the analogous
results when the correlation process is modelled by circular Brownian motion,
again under constant volatility and stochastic volatility, respectively. The
resulting correlation trajectories are broadly comparable to those obtained
under the von Mises specification, with the circular Brownian motion model
typically yielding slightly smoother paths, consistent with the absence of mean
reversion.

A comparison of Figures~\ref{fig:vmp_constvol} and \ref{fig:vmp_stochvol}
(von Mises correlation under constant and stochastic volatility,
respectively) shows that allowing for stochastic volatility does not
substantially alter the broad qualitative features of the estimated
correlation dynamics. In particular, the timing and direction of the major
correlation shocks associated with the onset of the COVID-19 pandemic are
preserved across the two specifications. This suggests that the dominant
features of correlation dynamics during periods of market stress are not an
artefact of assuming constant marginal volatilities. A similar pattern is observed when comparing Figures~\ref{fig:cbm_constvol}
and \ref{fig:cbm_stochvol}, corresponding to circular Brownian motion as the
correlation process. 

Overall, the empirical results suggest that the circular diffusion approach
yields interpretable time-varying correlation estimates for FX-equity pairs,
and that incorporating stochastic volatility does not alter the main qualitative
conclusions, while better accommodating periods of extreme market dislocation.

\section{Conclusion}\label{sec: conclusion}

In this article we studied diffusion processes on the unit circle with a focus on
statistical inference and applications to stochastic correlation modelling.
We considered two circular diffusions- the Brownian motion on the circle
and the von Mises process, a mean-reverting diffusion with the von Mises
distribution as its stationary law. A central contribution of the paper is the
derivation of an analytically tractable approximation to the transition density
of the von Mises process. To the best of our knowledge, this is the first
available analytical approximation to this transition density since the
original introduction of the von Mises diffusion by Kent
\cite{kent1975disc}, and it enables likelihood-based inference for
observed data.

Building on this approximation, we developed estimation procedures for circular
diffusions, studied their large-sample properties, and validated their
performance through simulation. We then proposed a novel stochastic correlation
model in which correlation is represented as the cosine of a latent angular
diffusion. This construction yields a bounded, interpretable correlation process
and allows both constant- and stochastic-volatility extensions. Estimation of
the resulting latent correlation paths was carried out via penalized likelihood,
with uncertainty quantified using a parametric bootstrap.

The proposed framework was illustrated through empirical applications to wind
direction data and to equity-foreign exchange market pairs during the COVID-19
pandemic. The results highlight the ability of circular diffusions to capture
time-varying correlation dynamics and to identify correlation shocks
during periods of market stress.

Several directions for future research arise naturally from this work. In
particular, the use of circular diffusions as correlation dynamics in
multi-asset derivative pricing models is a promising avenue. Extensions to higher-dimensional
settings and to dependence between correlation and asset dynamics also warrant
further investigation.

\section*
{Data Availability Statement}
The data that support the findings of this study are available as part of the R package \texttt{stochcorr}, available at \url{https://doi.org/10.32614/CRAN.package.stochcorr}.

\printbibliography

\appendix

\section{Proof of Lemma \ref{kents_construction}}

Let $X$ be a $d$-dimensional Riemannian manifold with metric tensor $g(x)$. Let $f(x)$ be a smooth, strictly positive density on $X$. The characteristic diffusion of $f$, i.e., the diffusion for which $f(x)$ is its stationary density is given by the pair $(G,B)$ (\cite{kent1978time}, Theorem 10.1):
\begin{equation}
	G=\Delta + \sum_i b_i(x)\frac{\partial}{\partial x_i}
\end{equation}
where $B$ is the periodic boundary condition. Here $G$ is the infinitesimal generator of the diffusion. More general forms of $B$ are also admissible, see p. 821 \cite{kent1978time}. Here $b_i(x)=\sum_{j=1}^dg^{-1}_{ij}\frac{\partial \log f}{\partial x_j}$ and $\Delta$ is the Laplace-Beltrami operator on $X$. The Laplace-Beltrami operator on $X$ is given by,

\begin{equation*}
\Delta =\delta^{1/2}(x)\sum_{i, j}\frac{\partial}{\partial x_i}\left[\delta^{-1/2}(x)g_{ij}(x)\frac{\partial}{\partial x_j}\right]
\end{equation*}
where $\delta(x)=\det g^{-1}(x)$. 

The Laplace-Beltrami operator on the unit circle is given by,
\begin{equation*}
\Delta=\frac{\partial^2}{\partial \theta^2}
\end{equation*}

Therefore the Kolmogorov backward equation is,

\begin{equation*}
    \frac{\partial p^*}{\partial t}=Gp^*
\end{equation*}

and hence,

\begin{equation*}
    \frac{\partial p^*}{\partial t}=\frac{\sigma^2}{2}\frac{\partial p^*}{\partial\theta^2}+\frac{\partial\log f}{\partial \theta}\frac{\partial p^*}{\partial\theta}
\end{equation*}

and collecting the coefficients for the drift and the diffusion term on the RHS, we obtain the SDE.

\section{Proof of Theorem \ref{thm:vmp}}

For a SDE,
\begin{equation}
	dY_t=\nu A(Y_t)dt+\sqrt{2\nu}dB_t
\end{equation}

the corresponding forward equation is given by,

\begin{equation}
	\frac{\partial f}{\partial \tau}=-\frac{\partial}{\partial y}[A(y)f]+\frac{\partial^2f}{\partial y^2}
\end{equation}

where $\tau$ is the nondimensional time, with $\tau=\nu t$. The stationary density  $f_e(y)\propto \exp \int A(z)dz$.

Then for  $f(\tau,y)=g(\tau,y)f_e(y)$, the adjoint equation is,

 \begin{equation}
	\frac{\partial g}{\partial \tau}=A(y)\frac{\partial g}{\partial y}+\frac{\partial^2 g}{\partial y^2}
\end{equation}

It follows that for $h=-\left(\frac{\partial }{\partial y}\right)\log g$, $h$ satisfies the PDE,
\begin{equation}
	\frac{\partial h}{\partial \tau}=\frac{\partial }{\partial y}\left\{A(y)h+\frac{\partial h}{\partial y}-h^2\right\}
\end{equation}

Consider an Ornstein-Uhlenbeck process,

\begin{equation*}
	dy_t=\nu \gamma(a-y_t)dt+\sqrt{2\nu} dB_t
\end{equation*}

For the OU process $A(y)=\gamma(a-y)$ and  $h$ is known in closed form for the OU process,

\begin{equation*}
	h(\tau,y)=\frac{\gamma\sqrt{q}(y-y_0)}{1-q}+\frac{\gamma\sqrt{q}(a-y)}{1+\sqrt{q}}
\end{equation*}

where $q = \exp(-2\gamma\tau)$. Notice that it can alternatively be expressed as,

\begin{equation}
    h(\tau,y) = \frac{\gamma\sqrt{q}(y - y_0)}{1 - q} + \frac{\sqrt{q}A(y)}{1 + \sqrt{q}}
\end{equation}
for $A(y)=\gamma(a-y)$. This form is used then as an ansatz for the solution to the forward equation for the von Mises process,

\begin{equation}\label{h_vm}
	h^*(\tau,y)=\frac{\gamma \sqrt{q}(y-y_0)}{1-q}+\frac{\sqrt{q}A(y)}{1+\sqrt{q}}+\sqrt{q}o(1)_{q=1}
\end{equation}

It is argued in \cite{martin2019analytical}, p. 5 that as $\tau\to 0$, i.e. $q\to 1$, the error term is of the order  $o(1)$, which is what $o(1)_{q=1}$ denotes. $\gamma$ is an artefact from the specification of the OU process. We need to identify a suitable value for $\gamma$ for the von Mises process. Towards this an expansion of the remainder term is considered,

\begin{equation*}
	h^*(\tau,y)=\frac{\gamma\sqrt{q}(y-y_0)}{1-q}+\frac{\sqrt{q}A(y)}{1+\sqrt{q}}+\frac{\sqrt{q}}{1+\sqrt{q}}\sum_{i=1}^{\infty}(1-\sqrt{q})^ib_i(y)
\end{equation*}

We note that,
\begin{equation}
b_1(y)=\frac{1}{\gamma(y-y_0)^2}\int_{y_0}^y(s-y_0)\left[\frac{d^2}{ds^2}(A(s)+\gamma s)+\left(A(s)+\frac{\gamma(s-y_0)}{2}\right)\left(\frac{d}{ds}(A(s)+\gamma s)\right)\right]ds
\end{equation}

Observe that $b_1(y)=0$ if  $A'(y)+\gamma=0$. This motivates the following choice for  $\gamma$,

 \begin{equation}
	 \gamma=\mathbb{E}_{f_{e}}[-A'(y)]
\end{equation}

where $f_e$ is the stationary distribution of the process. For the von Mises process,

\begin{equation}\label{vm_A}
	A(\theta)=-\frac{2\lambda}{\sigma^2}\sin(\theta-\mu)
\end{equation}

and $f_e=\frac{1}{2\pi I_0(\kappa)}\exp(\kappa\cos(\theta-\mu))$. Therefore,
\begin{equation}\label{vm_gamma}
	 \gamma=\int_{-\pi}^{\pi}\kappa\cos(\theta-\mu)\frac{\exp\left(\kappa\cos(\theta-\mu)\right)}{2\pi I_0\left(\kappa\right)}d\theta=\kappa\frac{I_1(\kappa)}{I_0(\kappa)}
\end{equation}

Plugging expressions (\ref{vm_A}) and (\ref{vm_gamma}) for von Mises process into (\ref{h_vm}) and we obtain,

\begin{equation}\label{unwrapped_tpd}
	f(\theta_{\tau};\theta_0)\propto \exp\left(\frac{-\gamma\sqrt{q}(\theta_{\tau}-\theta_0)^2}{2(1-q)}\right)\left(\frac{1}{2\pi I_0(\kappa)}\right)^{\frac{1-\sqrt{q}}{1+\sqrt{q}}}\exp\left(\frac{\kappa(\cos(\theta_{\tau}-\mu)-\sqrt{q}\cos(\theta_0-\mu))}{1+\sqrt{q}}\right)	
\end{equation}

To accommodate the periodic boundary conditions, we wrap (\ref{unwrapped_tpd}) to obtain the approximate transition density of the von Mises process,

\begin{align}
	p(\theta_t;\theta_0)&\propto\sum_{k=-\infty}^{\infty}\exp\left(\frac{-\gamma\sqrt{q}(\theta_t+2\pi k-\theta_0)^2}{2(1-q)}\right)\nonumber\\
			    &\left(\frac{1}{2\pi I_0(\kappa)}\right)^{\frac{1-\sqrt{q}}{1+\sqrt{q}}}\exp\left(\frac{\kappa(\cos(\theta_t-\mu)-\sqrt{q}\cos(\theta_0-\mu))}{1+\sqrt{q}}\right)	
\end{align}

where $\gamma=\kappa(I_1(\kappa)/I_0(\kappa))$ and $q=\exp(-\gamma\sigma^2t)$. We note that as $t\to \infty$,  $p(\theta_t;\theta_0)\to f_e$, which is the expected theoretical stationary distribution.

\section{Proof of Theorem \ref{thm:consistency}}

\begin{proof}Consider the von Mises process,
\begin{equation*}
d\theta_t = b(\theta_t;\lambda,\mu)\,dt + \sigma\,dW_t,
\qquad 
b(\theta;\lambda,\mu)=-\lambda\sin(\theta-\mu),
\end{equation*}
observed continuously on $[0,T]$.  
For fixed $\sigma>0$, the log-likelihood (up to an additive constant independent of $(\lambda,\mu)$) is
\begin{equation*}
L_T(\lambda,\mu)
=
\frac{1}{\sigma^2}\int_0^T b(\theta_t;\lambda,\mu)\,d\theta_t
-
\frac{1}{2\sigma^2}\int_0^T b(\theta_t;\lambda,\mu)^2\,dt .
\end{equation*}

Let $(\lambda_0,\mu_0)$ denote the true parameter value. The normalized log-likelihood may be written as
\begin{align*}
\frac{1}{T}L_T(\lambda,\mu)
&=
\frac{1}{\sigma^2}\frac{1}{T}\int_0^T 
b(\theta_t;\lambda,\mu)b(\theta_t;\lambda_0,\mu_0)dt
-\frac{1}{2\sigma^2}\frac{1}{T}\int_0^T b(\theta_t;\lambda,\mu)^2dt  \\
&\quad
+\frac{1}{\sigma}\frac{1}{T}
\int_0^T b(\theta_t;\lambda,\mu)dW_t .
\end{align*}

Define,
$$
M_T(\lambda,\mu)
=
\frac{1}{\sigma}
\int_0^T b(\theta_t;\lambda,\mu)\,dW_t .
$$

For each fixed $(\lambda,\mu)$, $M_T(\lambda,\mu)$ is a continuous martingale (\cite{ksendal2003}, Corollary 3.2.6)
with quadratic variation,
\begin{align*}
\langle M(\lambda,\mu)\rangle_T
&=\frac{1}{\sigma^2} \int_0^T b(\theta_t;\lambda,\mu)^2dt\\
&=\frac{1}{\sigma^2}\int_0^T \lambda^2\sin^2(\theta_t-\mu)dt
\end{align*}

Now consider,
\begin{align}
	&\sup_{(\lambda,\mu)\in K_{\lambda}\times K_{\mu}}\langle  M(\lambda,\mu) \rangle_T\nonumber\\
	=&\sup_{(\lambda,\mu)\in K_{\lambda}\times K_{\mu}} \frac{1}{\sigma^2}\int_0^T\lambda^2\sin^2(\theta_t-\mu)dt\nonumber\\
	\leq &\frac{1}{\sigma^2}\int_0^T\Lambda^2dt=\frac{\Lambda^2 T}{\sigma^2}\label{eq:quad_var}
\end{align}

where,
$$
\Lambda := \sup_{(\lambda,\mu)\in K_\lambda\times K_\mu} |\lambda|
< \infty,
$$

so that
$|b(\theta;\lambda,\mu)| \le \Lambda$
for all $\theta$ and all $(\lambda,\mu)\in K_\lambda\times K_\mu$.

From Ito's isometry and (\ref{eq:quad_var}),
\begin{equation}\label{eq:exp_quad_var}
\sup_{(\lambda,\mu)\in K_\lambda\times K_\mu}
\mathbb{E}\left[
\left|
\frac{1}{T}M_T(\lambda,\mu)
\right|^2
\right]=\sup_{(\lambda,\mu)\in K_\lambda\times K_\mu}
\mathbb{E}\left[
\frac{1}{T^2}\langle M_T(\lambda,\mu)\rangle_T\right]
\le \frac{\Lambda^2}{\sigma^2T},
\end{equation}

By Markov's inequality and (\ref{eq:exp_quad_var}),
$$
\mathbb{P}\left(\left|\frac{1}{T}M_T(\lambda,\mu)\right|>\varepsilon\right)\leq \frac{1}{\varepsilon^2}\mathbb{E}\left[\left|\frac{1}{T}M_T(\lambda,\mu)\right|^2\right]\leq \frac{\Lambda^2}{\varepsilon^2\sigma^2T}
$$

and hence,
\begin{equation*}
\sup_{(\lambda,\mu)\in K_\lambda\times K_\mu}
\left|
\frac{1}{T}M_T(\lambda,\mu)
\right|
\xrightarrow[T\to\infty]{\mathbb P}0.
\end{equation*}
Thus the stochastic integral term vanishes uniformly on compact subsets of
the parameter space.

We now consider the remaining deterministic integral terms.
By ergodicity of the von Mises process,
for any bounded measurable function $g$, (\cite{Kutoyants2004}, 1.55),
\begin{equation*}
\frac{1}{T}\int_0^T g(\theta_t)\,dt
\;\xrightarrow[T\to\infty]{\text{a.s.}}\;
\mathbb E_0[g(\Theta)],
\end{equation*}
where $\Theta$ has the invariant von Mises distribution with parameters corresponding to
$(\lambda_0,\mu_0)$.
Since the functions,
\begin{equation*}
g_1(\theta;\lambda,\mu)
=
b(\theta;\lambda,\mu)b(\theta;\lambda_0,\mu_0),
\qquad
g_2(\theta;\lambda,\mu)
=
b(\theta;\lambda,\mu)^2
\end{equation*}
are bounded and continuous jointly in $(\theta,\lambda,\mu)$,
the ergodic theorem applies.
Therefore,
\begin{equation}\label{eq:likelihood_conv}
\frac{1}{T}L_T(\lambda,\mu)
\xrightarrow[T\to\infty]{\mathbb P}
 Q(\lambda,\mu),
\end{equation}
uniformly on compact subsets, where
\begin{equation*}
 Q(\lambda,\mu)
=
\frac{1}{\sigma^2}
\mathbb E_0\!\big[
b(\Theta;\lambda,\mu)b(\Theta;\lambda_0,\mu_0)
\big]
-
\frac{1}{2\sigma^2}
\mathbb E_0\!\big[
b(\Theta;\lambda,\mu)^2
\big].
\end{equation*}

Expanding the square,
\begin{equation*}
Q(\lambda,\mu)
=
\text{\rm const}
-
\frac{1}{2\sigma^2}
\mathbb E_0\!\Big[
\big(
b(\Theta;\lambda,\mu)-b(\Theta;\lambda_0,\mu_0)
\big)^2
\Big].
\end{equation*}

Hence $Q(\lambda,\mu)$ is uniquely maximized at
$(\lambda_0,\mu_0)$.
If,
\begin{equation*}
b(\theta;\lambda,\mu)
=
b(\theta;\lambda_0,\mu_0)
\quad\forall\theta,
\end{equation*}
then,
\begin{equation*}
\lambda\sin(\theta-\mu)
=
\lambda_0\sin(\theta-\mu_0)
\quad\forall\theta.
\end{equation*}
This implies,
\begin{equation*}
(\lambda\cos\mu-\lambda_0\cos\mu_0)\sin\theta
-
(\lambda\sin\mu-\lambda_0\sin\mu_0)\cos\theta
=0
\quad\forall\theta.
\end{equation*}
By linear independence of $\sin\theta$ and $\cos\theta$, we obtain,
\begin{equation*}
\lambda\cos\mu=\lambda_0\cos\mu_0,
\qquad
\lambda\sin\mu=\lambda_0\sin\mu_0.
\end{equation*}
Since $\lambda,\lambda_0>0$, this yields $\lambda=\lambda_0$ and
$\mu=\mu_0 \pmod{2\pi}$.

The following also follows from \cite{Vaart1998}, Theorem 5.7, however we state the argument here for completeness. Fix $\varepsilon>0$ and set
\begin{equation*}
K := K_\lambda \times K_\mu,
\qquad
K_\varepsilon
:=
\bigl\{(\lambda',\mu')\in K :
|(\lambda',\mu')-(\lambda,\mu)| \ge \varepsilon\bigr\}.
\end{equation*}
Since $K$ is compact, the set $K_\varepsilon$ is compact as a closed subset
of $K$.
Moreover, the function $ Q$ is continuous on $K$ and admits a unique
maximizer at $(\lambda,\mu)$.
Hence $Q$ attains its maximum on $K_\varepsilon$ at some point
$(\lambda_\varepsilon,\mu_\varepsilon)\in K_\varepsilon$, and,
\begin{equation*}
Q(\lambda_\varepsilon,\mu_\varepsilon)
< Q(\lambda,\mu).
\end{equation*}

For every $\varepsilon>0$ there exists
$\eta(\varepsilon)>0$ such that,
\begin{equation*} Q(\lambda,\mu)
=
\sup_{(\lambda',\mu')\in K_\varepsilon}Q(\lambda',\mu')
+\eta(\varepsilon).
\end{equation*}

Moreover from (\ref{eq:likelihood_conv}),
\begin{equation*}
\sup_{(\lambda',\mu')\in K}
\bigl|
T^{-1}L_T(\lambda',\mu')-Q(\lambda',\mu')
\bigr|
\xrightarrow[T\to\infty]{\mathbb P}0.
\end{equation*}
Fix $\delta\in(0,\eta(\varepsilon)/3)$ and consider the event,
\begin{equation*}
V_T
=
\left\{
\sup_{(\lambda',\mu')\in K}
\bigl|
T^{-1}L_T(\lambda',\mu')-Q(\lambda',\mu')
\bigr|
\le \delta
\right\}.
\end{equation*}

On $V_T$, for all $(\lambda',\mu')\in K$,
\begin{equation}\label{eq:lik_bounds} 
	Q(\lambda',\mu')-\delta
\le
T^{-1}L_T(\lambda',\mu')
\le Q(\lambda',\mu')+\delta.
\end{equation}

In particular,
\begin{equation*}
T^{-1}L_T(\lambda,\mu)
\ge Q(\lambda,\mu)-\delta
=
\sup_{(\lambda',\mu')\in K_\varepsilon} Q(\lambda',\mu')
+\eta(\varepsilon)-\delta.
\end{equation*}

Since $\delta<\eta(\varepsilon)/3$, we have
$\eta(\varepsilon)-\delta>\delta$, and therefore
\begin{equation}\label{eq:Q_bounds}
	Q(\lambda,\mu)-\delta
>
\sup_{(\lambda',\mu')\in K_\varepsilon} Q(\lambda',\mu')
+\delta.
\end{equation}

From (\ref{eq:lik_bounds}) and (\ref{eq:Q_bounds}), it follows that on $V_T$,
\begin{equation}\label{eq:L_bound}
T^{-1}L_T(\lambda,\mu)
>
\sup_{(\lambda',\mu')\in K_\varepsilon}
T^{-1}L_T(\lambda',\mu').
\end{equation}

Suppose, for contradiction, that there exists a maximizer
$(\hat\lambda,\hat\mu)$ of $L_T$ over $K$ such that
$|(\hat\lambda,\hat\mu)-(\lambda,\mu)|\ge\varepsilon$.
Then $(\hat\lambda,\hat\mu)\in K_\varepsilon$, and by definition of a
maximizer,
\begin{equation*}
T^{-1}L_T(\hat\lambda,\hat\mu)
=
\sup_{(\lambda',\mu')\in K}
T^{-1}L_T(\lambda',\mu')
\ge
T^{-1}L_T(\lambda,\mu),
\end{equation*}
which contradicts (\ref{eq:L_bound}).
Hence, on $V_T$, every maximizer $(\hat\lambda,\hat\mu)$ of
$L_T$ over $K$ satisfies
\begin{equation*}
|(\hat\lambda,\hat\mu)-(\lambda,\mu)|<\varepsilon.
\end{equation*}

Therefore,
\begin{equation*}
\mathbb P_{\lambda,\mu}
\!\left(
|(\hat\lambda,\hat\mu)-(\lambda,\mu)|\ge\varepsilon
\right)
\le
\mathbb P(V_T^{\,c}),
\qquad
\forall(\lambda,\mu)\in K.
\end{equation*}
Since,
\begin{equation*}
\sup_{(\lambda,\mu)\in K}
\mathbb P(V_T^{\,c})
\to 0\text{ as }T\to\infty,
\end{equation*}
we conclude that
\begin{equation*}
\sup_{(\lambda,\mu)\in K}
\mathbb P_{\lambda,\mu}
\!\left(
|(\hat\lambda,\hat\mu)-(\lambda,\mu)|\ge\varepsilon
\right)
\to 0,
\end{equation*}
which establishes uniform consistency on $K$.

Finally, under continuous observation,
\begin{equation*}
[\theta]_T=\sigma^2T \quad \text{a.s.},
\end{equation*}
so the maximum likelihood estimator
$\hat\sigma^2=[\theta]_T/T$ equals $\sigma^2$ almost surely for all $T$,
yielding uniform consistency of $\hat\sigma$.
\end{proof}

\section{Proof of Theorem \ref{thm:normality}}

\begin{proof}
$$
d\theta_t=-\lambda\sin(\theta_t-\mu)dt+\sigma dW_t
$$

is an ergodic diffusion on the compact state space $\mathbb S^1$ with smooth
drift and constant diffusion coefficient $\sigma>0$. Writing it in the form
of \cite{Kutoyants2004},
$$
dX_t=S(\beta,X_t)dt+a(X_t)dW_t,
\qquad
\beta=(\lambda,\mu),
$$

we have $X_t=\theta_t$, $a(x)\equiv\sigma$, and
$S(\beta,x)=-\lambda\sin(x-\mu)$.
The function $S(\beta,\cdot)$ is smooth and bounded uniformly on compact
subsets of the parameter space, and the diffusion coefficient is bounded away
from zero. Moreover, identifiability holds once $\lambda>0$ and
$\mu\in[0,2\pi)$ is fixed modulo $2\pi$. Hence the regularity and ergodicity
conditions $A0(\Theta)$ and $A$ of \cite{Kutoyants2004} are satisfied
on any compact subset of $(0,\infty)\times[0,2\pi)$.

Fix $(\lambda_0,\mu_0)$ and consider local perturbations
$\lambda=\lambda_0+u_1/\sqrt{T}$ and $\mu=\mu_0+u_2/\sqrt{T}$.
By Lemma~2.9 of \cite{Kutoyants2004}, the likelihood ratio admits a uniformly
locally asymptotically normal expansion,
\[
L_T(\lambda,\mu)-L_T(\lambda_0,\mu_0)
=
u^\top\Delta_T
-\frac12 u^\top I(\lambda_0,\mu_0)u
+o_{P_0}(1),
\]
where the central sequence is given by
\[
\Delta_T
=
\frac{1}{\sqrt{T}}\int_0^T
\frac{\dot S(\lambda_0,\mu_0;\theta_t)}{\sigma}\,dW_t
=
\frac{1}{\sigma\sqrt{T}}
\int_0^T
\begin{pmatrix}
-\sin(\theta_t-\mu_0)\\[2pt]
\lambda_0\cos(\theta_t-\mu_0)
\end{pmatrix}
\,dW_t.
\]
By the martingale central limit theorem,
$\Delta_T$ converges in distribution to a centered normal random vector with
covariance matrix $I(\lambda_0,\mu_0)$.

For ergodic diffusions of this form, the Fisher information matrix is,
\[
I(\lambda_0,\mu_0)
=
\mathbb E_{\lambda_0,\mu_0}
\!\left[
\frac{\dot S(\lambda_0,\mu_0;\Theta)\dot S(\lambda_0,\mu_0;\Theta)^\top}
{\sigma^2}
\right],
\]
where $\Theta$ is distributed according to the invariant law of the diffusion.
In the present case this invariant law is the von Mises distribution with parameters $\mu_0$ and $\kappa=2\lambda_0/\sigma^2$.
By symmetry, the cross term vanishes, and the information matrix is diagonal.
Using the definition of modified Bessel function of the first kind we obtain,
$$
\mathbb E[\sin^2(\Theta-\mu_0)]
=
\frac{I_0(\kappa)-I_2(\kappa)}{2I_0(\kappa)},
\qquad
\mathbb E[\cos^2(\Theta-\mu_0)]
=
\frac{I_0(\kappa)+I_2(\kappa)}{2I_0(\kappa)}.
$$

Consequently,
\[
I_{\lambda\lambda}
=
\frac{I_0(\kappa)-I_2(\kappa)}{2\sigma^2 I_0(\kappa)},
\qquad
I_{\mu\mu}
=
\frac{\lambda_0^2\bigl(I_0(\kappa)+I_2(\kappa)\bigr)}
{2\sigma^2 I_0(\kappa)}.
\]

Since the model is uniformly locally asymptotically normal and the regularity
conditions of \cite{Kutoyants2004} hold, Theorem~2.8 therein applies. It follows
that the maximum likelihood estimator is uniformly asymptotically normal on
compact subsets of the parameter space, with asymptotic covariance matrix equal
to the inverse of the Fisher information. Inverting the diagonal matrix
$I(\lambda_0,\mu_0)$ yields,
$$
\sqrt{T}(\hat{\lambda}-\lambda_0)
\xrightarrow{d}
N\left(0,\frac{2\sigma^2 I_0(\kappa)}{I_0(\kappa)-I_2(\kappa)}\right),
\qquad
\sqrt{T}(\hat{\mu}-\mu_0)
\xrightarrow{d}
N\left(0,\frac{2\sigma^2 I_0(\kappa)}
{\lambda_0^2\bigl(I_0(\kappa)+I_2(\kappa)\bigr)}\right),
$$

uniformly on compact subsets. 
\end{proof}

\section{Proof of Theorem \ref{thm:strong}}

\subsection{Preliminaries}

In this section we consider the existence of strong solution of the following system of SDEs,

\begin{align}
	dS_t^{(1)}&=\mu^{(1)}S_t^{(1)}dt+\sigma^{(1)}S_t^{(1)}dW_t^{(1)}\label{sde1}\\
	dS_t^{(2)}&=\mu^{(2)}S_t^{(2)}dt+\sigma^{(2)}S_t^{(2)}(\rho_t dW_t^{(1)}+\sqrt{1-\rho_t^2}dW_t^{(2)})\label{sde2}\\
	d\theta_t&=-\lambda\sin(\theta_t-\mu)dt+\sigma dW_t^{(3)}\label{sde3}
\end{align}

Here $W_t^{(1)},W_t^{(2)},W_t^{(3)}$ are mutually independent Brownian motions and $\rho_t=\cos\theta_t$. The existence of strong solution of (\ref{sde1}) and (\ref{sde3}) follows from the standard theory of Lipschitz coefficients. We will prove that (\ref{sde2}) admits a strong solution. We will require the following version of the Gr{\"o}nwall's inequality,
\begin{lemma}\label{lemma: gronwall}
	Let $y_n(t)$ be a sequence of nonnegative functions such that for some constants  $B,C<\infty$,  $y_0(t)\leq C,\forall t\leq T$ and  $y_{n+1}(t)\leq B\int_0^{t}y_n(s)ds<\infty,\forall t\leq T, n\in \mathbb{N}$. Then,
	$$
	y_n(t)\leq \frac{CB^nt^n}{n!}, \forall t\leq T
	$$
	\end{lemma}
\subsection{Strong solution}

We will rewrite (\ref{sde2}) in the following manner to ease on the notation,
\begin{equation}\label{sde}
	dS(t)=\mu S(t)dt +\sigma S(t)(\rho(t)dW^{(1)}(t)+\sqrt{1-\rho^2(t)}dW^{(2)}(t))
\end{equation}
Consider the initial condition $S(0)=S$ and a compact time interval $[0,T]$. We next consider the Picard's iterations of (\ref{sde}). Therefore for some $t\in [0,T]$,
\begin{align*}
	S_0(t)&\equiv S\\
	S_{n+1}(t)&=S+\int_0^t\mu S_n(u)du+\int_0^t\sigma S_n(u)(\rho(u)dW^{(1)}(u)+\sqrt{1-\rho^2(u)}dW^{(2)}(u)), \forall n\in  \mathbb{N}
\end{align*}

We will next show that $S_n(t)$ converges in  $L^2$.
\begin{align*}
	&\mathbb{E}[(S_{n}(t)-S_{n-1}(t))^2]\\
	=&\mathbb{E}\left[\left(\int_0^t\mu (S_{n}(u)-S_{n-1}(u))du+\int_0^t\sigma (S_{n}(u)-S_{n-1}(u))(\rho(u)dW^{(1)}(u)+\sqrt{1-\rho^2(u)}dW^{(2)}(u)\right)^2\right]\\
	\leq&4\mathbb{E}\left[\left(\int_0^t\mu(S_{n}(u)-S_{n-1}(u))du\right)^2\right]+4\mathbb{E}\left[\left(\int_0^t\sigma\rho(u)(S_{n}(u)-S_{n-1}(u))dW^{(1)}(u)\right)^2\right]\\
	    &+4\mathbb{E}\left[\left(\int_0^t\sigma\sqrt{1-\rho^2(u)}(S_{n}(u)-S_{n-1}(u))dW^{(1)}(u)\right)^2\right]\\
	    &\text{(Follows from $(x+y+z)^2\leq 4x^2+4y^2+4z^2$)}\\
	\leq&4\mathbb{E}\left[t\int_0^t\mu^2(S_{n}(u)-S_{n-1}(u))^2du\right]+4\mathbb{E}\left[\int_0^t\sigma^2\rho^2(u)(S_{n}(u)-S_{n-1}(u))^2du\right]\\
	    &+4\mathbb{E}\left[\int_0^t\sigma^2(1-\rho^2(u))(S_{n}(u)-S_{n-1}(u))^2 du\right]\\
	    &\text{(Applying Cauchy-Schwarz inequality to the first term and}\\
	    &\text{It{\^o}'s isometry to the second and third term)}\\
	\leq &4K(t+1)\mathbb{E}\left[\int_0^t(S_{n}(u)-S_{n-1}(u))^2du\right]\\
		 &\text{($\exists K>\max(\mu^2,\sigma^2)$)}\\
	=&4K(t+1)\int_0^t\mathbb{E}\left[(S_{n}(u)-S_{n-1}(u))^2\right]du
\end{align*}

From Lemma \ref{lemma: gronwall}, 
	\begin{equation}
		\mathbb{E}\left[(S_{n+1}(u)-S_n(u))^2\right]\leq \frac{C(4K(t+1))^n}{n!}		
	\end{equation}
Therefore, $\exists S(t)$ such that,
\begin{equation}
 \lim_{n\to \infty}\mathbb{E}[(S_{n}(t)-S(t))^2]=0
\end{equation}
for $t\in [0,T]$.

It also follows that,
\begin{align}
	\lim_{n\to \infty}&\mathbb{E}\left[(\mu(S_n(t)-S(t)))^2\right]=0\label{term1}\\
	\lim_{n\to \infty}&\mathbb{E}\left[(\sigma\rho(u)(S_n(t)-S(t)))^2\right]\leq\lim_{n\to \infty}\mathbb{E}\left[(\sigma(S_n(t)-S(t)))^2\right]=0\label{term2}\\
	\lim_{n\to \infty}&\mathbb{E}\left[(\sigma\sqrt{1-\rho^2(u)}(S_n(t)-S(t)))^2\right]\leq\lim_{n\to \infty}\mathbb{E}\left[(\sigma(S_n(t)-S(t)))^2\right]=0\label{term3}\\
			  &\text{((\ref{term2}) and (\ref{term3}) follows since $\abs{\rho(u)}\leq 1$)}\nonumber
\end{align}
We will next show that $S(t)$ satisfies (\ref{sde}). 
Consider,
\begin{align}
	&\lim_{n\to\infty}\mathbb{E}\left[\left(\int_0^t \mu(S_n(u)-S(u))du\right)^2\right]\nonumber\\
	\leq&\lim_{n\to \infty}\mathbb{E}\left[t\int_0^t(\mu(S_n(u)-S(u)))^2du\right]\nonumber\\
	    &\text{(Using Cauchy-Schwarz inequality)}\nonumber\\
	=&\lim_{n\to \infty}\int_0^tt\mathbb{E}\left[(\mu(S_n(u)-S(u)))^2\right]du=0\label{drift}\\
	 &\text{(From (\ref{term1}))}\nonumber
\end{align}
From (\ref{drift}) it follows that the drift term of (\ref{sde}) converges in $L^2$.
Next consider,
\begin{align}
	&\lim_{n\to\infty}\mathbb{E}\left[\left(\int_0^t \sigma\rho(u)(S_n(u)-S(u))dW^{(1)}(u)\right)^2\right]\nonumber\\
	=&\lim_{n\to \infty}\mathbb{E}\left[\int_0^t(\sigma\rho(u)(S_n(u)-S(u)))^2du\right]\nonumber\\
	    &\text{(Using It{\^o}'s isometry)}\nonumber\\
	\leq&\lim_{n\to \infty}\int_0^t\mathbb{E}\left[(\sigma(S_n(u)-S(u)))^2\right]du=0\label{diff1}\\
	 &\text{(From (\ref{term2}) and $\abs{\rho(u)}\leq 1$)}\nonumber
\end{align}
Similarly, we can show that,
\begin{equation}\label{diff2}
	\lim_{n\to\infty}\mathbb{E}\left[\left(\int_0^t \sigma\sqrt{1-\rho^2(u)}(S_n(u)-S(u))dW^{(2)}(u)\right)^2\right]=0
\end{equation}

From (\ref{diff1}) and (\ref{diff2}) it follows that the diffusion terms of (\ref{sde}) converge in $L^2$.
\end{document}